\documentclass{article}
\usepackage{amssymb}
\usepackage{rotating}
\usepackage{url}

\newtheorem{termenadenn}{ \ \ D\'efinition}[section]
\newtheorem{kinnig}[termenadenn]{ \ \ Proposition}
\newtheorem{teorem}[termenadenn]{ \ \ Th\'eor\`eme}
\newtheorem{merkadenn}[termenadenn]{ \ \ Remarque}
\newtheorem{lemma}[termenadenn]{ \ \ Lemme}
\newtheorem{heuliadenn}[termenadenn]{ \ \ Corollaire}

\newtheorem{termenadenn_lechel}{ \ \ D\'efinition}[subsection]
\newtheorem{kinnig_lechel}[termenadenn_lechel]{ \ \ Proposition}

\newtheorem{merkadenn_lechel}[termenadenn_lechel]{ \ \ Remarque}

\newcounter{biblio}

\newcommand{\hed}{\hspace*{.7cm}}
\newcommand{\supp}{{\sf supp}}
\newcommand{\dom}{{\sf dom}}

\begin{document}

\title{Invariance combinatoire des
polyn\^omes de Kazhdan-Lusztig sur les intervalles
partant de l'origine}

\author{ Ewan Delanoy 
\footnote{Institut Girard D\'esargues \newline
UMR 5028 CNRS \newline
Universit\'e Lyon 1 \newline
 69622 Villeurbanne Cedex France \newline
delanoy@igd.univ-lyon1.fr
}
}

\maketitle

\begin{abstract}

{ We show that for Bruhat intervals starting from
the origin in Coxeter groups the
conjecture of Lusztig and Dyer holds, that is, the $R$-polynomials
and the Kazhdan-Lusztig polynomials defined on 
$[e,u]$ only depend on 
the isomorphism type of $[e,u]$. To achieve this we
use the purely poset-theoretic notion of special matching.
Our approach is essentially a synthesis of the
explicit formula for special matchings discovered by Brenti and
the general special matching machinery developed by Du Cloux. }
\end{abstract}

\bigskip
 
\section{Introduction}\label{intro}

\vspace*{2cm}

 \hed \`A  l'heure actuelle la question pos\'ee 
ind\'ependamment par
Dyer [\ref{Dyer_1987}] et Lusztig  de savoir si
le polyn\^ome de Kazhdan-Lusztig $P_{u,v}$ d\'efini sur
un intervalle de Bruhat $[u,v]$ ne d\'epend en fait
que de la classe d'isomorphisme du poset $[u,v]$,
reste un probl\`eme ouvert dans le cas g\'en\'eral 
(voir [\ref{Humphreys}] pour toutes les
d\'efinitions concernant l'ordre de Bruhat et
les polyn\^omes $P$ et $R$ associ\'es \`a un groupe
de Coxeter ). On peut reformuler le probl\`eme
de la mani\`ere suivante : est-ce que pour tout
isomorphisme de posets $\psi$ entre deux intervalles
de Bruhat, $\psi:{[u,v]}_W \to {[u',v']}_{W'}$, $\psi$
pr\'eserve les polyn\^omes de Kazhdan-Lusztig i.e. \newline

\hspace*{3.5cm} {$\forall x,y \in [u,v], \
P_{x,y}=P_{\psi(x),\psi(y)} $} \hfill (1.1) \newline

Brenti [\ref{Brenti_1994}]
a trait\'e le cas adih\'edral (i.e. le cas o\`u le
poset $[u,v]$ ne contient pas d'intervalle isomorphe
au poset ${\mathfrak S}_3$ vu comme groupe de Coxeter muni
de l'ordre de Bruhat). Nous montrons
dans cet article que le cas particulier cor\-res\-pondant
\`a $u=u'=e$ est vrai : \newline

\hspace*{3.5cm} {$\forall x,y \in [e,v], \
P_{x,y}=P_{\psi(x),\psi(y)} $} \hfill (1.2) \newline

 Des sous-cas de ce cas particulier ont d\'eja \'et\'e
d\'emontr\'es : la r\'ef\'erence [\ref{du_Cloux_2003}] 
traite le cas o\`u $W$ est tel que
 toutes les composantes connexes de son graphe de Coxeter
sont des arbres ou de type 
${\widetilde{A}}_n$
et [\ref{Brenti_2004}] traite le cas o\`u $W$ et $W'$
sont de type $A_n$.
On sait qu'il existe un algorithme
permettant de calculer les polyn\^omes $P$
\`a partir d'une famille de polyn\^omes plus 
\'el\'ementaires, les polyn\^omes $R$ (cf. par
exemple \ref{Brenti_2004}, th\'eor\`eme 2.6.iv)) ; cet
algorithme qui n'utilise que la structure de poset
des intervalles de Bruhat montre que (1.2) \'equivaut 
\`a \newline

\hspace*{3.5cm} {$\forall x,y \in [e,v], \
R_{x,y}=R_{\psi(x),\psi(y)} $} \hfill (1.3) \newline

 Pour montrer un r\'esultat d'invariance
par isomorphisme de posets, on cherche tout naturellement
\`a d\'efinir de mani\`ere purement combinatoire les
polyn\^omes $R$ ; ce qui a conduit
Du Cloux et Brenti \`a la notion de ``couplage distingu\'e'',
que nous explicitons un peu plus loin. \'Etant donn\'e un
g\'en\'erateur $s$, les op\'erateurs de multiplication \`a gauche
et \`a droite par $s$ constituent des exemples
fondamentaux de couplages distingu\'es ; nous
les appelons des couplages distingu\'es {\it de multiplication}. 
  Par un raisonnement commun \`a [\ref{du_Cloux_2000}, d\'efinition 6.5] et 
[\ref{Brenti_2004}, corollaire 5.3], on montre que
l'assertion (1.3)
est impliqu\'ee par l'assertion suivante 
portant les couplages distingu\'es $\phi$
d'un intervalle partant de l'origine $[e,v]$ 
(bien connue lorsque $\phi$ est un couplage de multiplication):

\begin{displaymath}
\begin{array}{c}
\forall x,y \in [e,v], {\rm \ tels \ que \ }
x\lhd \phi(x), y \lhd \phi(y), \\
\hspace*{3cm}
\left\lbrace
\begin{array}{l}
R_{\phi (x), \phi (y)}=R_{x,y} \\
R_{x, \phi (y)}=(q-1)R_{x,y}+qR_{\phi(x),y}
\end{array}\right.\hspace*{2.5cm}(1.4)
\end{array}
\end{displaymath} 

 C'est ce r\'esultat que nous d\'emontrons dans
cet article (corollaire \ref{final_result}). \newline

 L'id\'ee de d\'epart de la d\'emonstration, d\'eja contenue
dans [\ref{Brenti_2004}], est la suivante~: prenons
$(x,y)$ comme en (1.4), et
$s\in S$ dans l'ensemble de descente \`a gauche de $y$ tel que
$\phi$ commute avec avec la multiplication \`a gauche par $s$.
Alors la formule (1.4) pour $(x,y)$ se d\'eduit de formules (1.4)
correspondant \`a des $(x',y')$ avec $l(y')<l(y)$, ce qui permet
de raisonner par r\'ecurrence sur la longueur de $y$ (proposition
\ref{descent_formulas}).  Il n'est pas vrai en g\'en\'eral que
tout $y \in W$ admette une telle r\'eduction (\`a gauche ou \`a
droite). Mais ce sera vrai pour tout \'el\'ement ``suffisamment
grand''. Plus pr\'ecisement, disons que $y$ est {\it plein} si
$[e,y]$ contient tous les \'el\'ements dih\'edraux de $W$
(d\'efinition \ref{full_elements}). Alors nous d\'emontrons
{\it in fine } que si $W$ n'est pas dih\'edral, tout
\'el\'ement plein de $W$ est r\'eductible au sens ci-dessus. 
Si $w$ est non plein, il est prouv\'e dans [\ref{du_Cloux_2003}]
 que l'intervalle $[e,w]$ est isomorphe \`a un intervalle
$[e,w']$ dans un autre groupe de Coxeter $W'$ 
``plus petit'' en un sens convenable par un isomorphisme
pr\'eservant les polyn\^omes $R$, avec $w'$ plein, ce qui permet
de faire une r\'ecurrence sur la taille du groupe de Coxeter.\newline

 Dans tous les raisonnements , les \'el\'ements dih\'edraux 
du groupe jouent un r\^ole essentiel. \`A la section \ref{resgen}
nous montrons qu'un couplage est enti\`erement caract\'eris\'e
par sa restriction aux \'el\'ements dih\'edraux
de son domaine de d\'efinition,
et m\^eme par sa restriction \`a l'ensemble $P$ des \'el\'ements
dih\'edraux principaux (th\'eor\`eme \ref{mw_is_a_cartesian_product}).
R\'eciproquement tout couplage $\phi$ d\'efini sur $P$ se prolonge
de mani\`ere unique en un couplage 
(encore not\'e $\phi$) dont le domaine est maximal~; ce domaine uniquement
d\'efini est not\'e $\dom(\phi)$. De m\^eme, 
\`a la section \ref{big_enough} on verra que la commutation
d'un couplage avec l'op\'erateur de multiplication
par un g\'en\'erateur se lit sur $P$ (proposition
\ref{commutativity_on_ps_suffices} ). \newline

 Pour que  $\dom(\phi)$ contienne un \'el\'ement plein,
il faut que la restriction de $\phi$ \`a chaque sous-groupe
dih\'edral principal ne soit ``pas trop \'eloign\'ee'' d'un couplage
de multiplication : nous verrons \`a la section \ref{gen}
qu'il existe au plus un sous-groupe dih\'edral principal $D$
tel que la restriction de $\phi$ \`a $D$ ne soit pas un couplage
de multiplication, et que m\^eme sur ce $D$ $\phi$ doit encore
partager avec les couplages de multiplication
certaines conditions de r\'egularit\'e. \newline

  Sur le plan technique, une id\'ee essentielle 
consiste \`a mettre en \'evidence des ``obstructions'' 
(\'el\'ement minimaux du compl\'e\-mentaire de $\dom(\phi)$)
chaque fois que $\phi$ n'est pas un couplage de multiplication.
Par exemple, si $a=\phi(e)$ et $x_0 \in P$ est un \'el\'ement 
minimal tel que $\phi(x_0) \neq x_0a$, on peut
exhiber des obstructions obtenues en ins\'erant un caract\`ere
bien choisi dans une \'ecriture r\'eduite de $x_0$
(ceci est illustr\'e par les propositions \ref{core_of_nd}
et \ref{obs1_d}). Comme le domaine $\dom(\phi)$ est filtrant d\'ecroissant,
chaque nouvelle obstruction impose une diminution cons\'equente
de $\dom(\phi)$, pour finalement emp\^echer $\dom(\phi)$
de contenir un \'el\'ement plein lorsque $\phi$ est trop
diff\'erent d'un couplage de multiplication, ce qui permet
de faire aboutir le raisonnement.\newline   

 Il est remarquable que les seules obstructions dont 
nous ayons besoin pro\-viennent toutes de sous-groupes de
rang 3 de $W$. \`A la section
\ref{not_too_big} nous d\'ecrivons les obstructions en
quelque sorte les plus simples que
l'on puisse rencontrer et la r\'eduction correspondante du
domaine $Q$, qui apparait dans le cas dit ``crois\'e'',
qui r\`egle d\'eja le cas des groupes de Coxeter simplement
enlac\'es (cf corollaire \ref{simply_laced_case}). La section
\ref{rank_three}, consacr\'ee au cas du rang 3, fournit des
armes pour l'attaque du cas g\'en\'eral (section \ref{gen}). 
 La d\'ecouverte des obstructions en rang 3
a \'et\'e largement guid\'ee par des calculs \'effectu\'es \`a
l'aide d'une version sp\'ecialis\'ee du programme
{\tt Coxeter} [\ref{du_Cloux_Coxeter}].\newline

 Plan de la suite de l'article :
\begin{trivlist}
\item[\S\ref{resgen}.] R\'esultats g\'en\'eraux
\item[\S\ref{exit_r}.] Formules de descente pour
les polyn\^omes {\it R}
\item[\S\ref{big_enough}.] Crit\`eres de r\'egularit\'e
\item[\S\ref{not_too_big}.] R\'eduction du domaine dans le cas crois\'e
\item[\S\ref{rank_three}.] \'Etude partielle du cas des groupes de rang 3
\item[\S\ref{gen}.] Cas g\'en\'eral
\end{trivlist}

\bigskip

\section{R\'esultats g\'en\'eraux}\label{resgen}

 \hed  Soit $(P,<)$ un poset. On \'ecrit $x \triangleleft y$
pour exprimer que $x<y$ et qu'il n'existe pas de $z$ tel
que $x<z<y$. On dit alors que
$x$ est un {\bf coatome} de $y$ ;
on note {\bf coat(y)} l'ensemble des coatomes d'un
\'el\'ement $y \in P$. Tous les 
posets consid\'er\'es ici 
seront
{\bf gradu\'es}, i.e. munis d'une fonction $l \ : \ 
P \rightarrow \Bbb{N}$ v\'erifiant $l(y)=l(x)+1$
d\`es que $x \triangleleft y$ (en fait, ne sont 
consid\'er\'es ici que des posets gradu\'es en g\'eneral
sans aucune autre restriction ou bien, dans le cadre
d'un syst\`eme de Coxeter $(W,S)$, on regardera toujours
$W$ comme \'etant muni de sa structure de poset gradu\'e
provenant de l'ordre de Bruhat et de la fonction longueur
usuelle). 

 Soit $\phi \ : \ P \rightarrow P$ une application. Nous 
dirons que $\phi$ est un {\bf couplage distingu\'e}
si les trois conditions suivantes sont v\'erifi\'ees:

\begin{trivlist}
\item[(i)]
$\phi$ est involutive ( $ \forall u \in P, \ \phi(\phi(u))=u $   )
\item[(ii)]
$\forall u \in P, \ u \triangleleft \phi(u) \ {\rm ou} \ 
\phi(u) \triangleleft u$
\item[(iii)]
$
\forall u \in P, \ 
(u \triangleleft \phi(u)) \Rightarrow
( \ coat(\phi(u))= \lbrace u  \rbrace \cup \lbrace \phi(v) \ ; \ 
v \triangleleft u, \
v\triangleleft  \phi(v) \rbrace \  ).
$
\end{trivlist}

 La condition $(iii)$ est la plus significative, les autres
ne font que poser le cadre. Nous utiliserons l'abbr\'evation
\begin{displaymath}
Z( \phi, u)=\lbrace u  \rbrace \cup \lbrace \phi(v) \ ; \ 
v \triangleleft u, \
v\triangleleft  \phi(v) \rbrace
\end{displaymath}

 Si la terminologie est due \`a Brenti [\ref{Brenti_2004}], 
le choix de la d\'efinition
(parmi un certain nombre qui sont \'equivalentes) vient plut\^ot de
du Cloux [\ref{du_Cloux_2000}]. Les parties $(i)$ et $(ii)$ de la
d\'efinitions sont communes \`a [\ref{Brenti_2004}] et \`a 
[\ref{du_Cloux_2000}] ; la
partie $(iii)$ par contre n'est \'enonc\'ee explicitement
ni dans [\ref{Brenti_2004}] ni dans [\ref{du_Cloux_2000}], 
mais il est facile de voir
qu'elle est \'equivalente aux versions donn\'ees dans
chacun de ces articles.

 \hspace{1cm} Si l'on prend $Q$, une 
partie {\bf filtrante d\'ecroissante } de $P$
(i.e. v\'erifiant  
$\forall (q,x) \in Q \times P, \ 
(x \leq q) \Rightarrow (x \in Q)
$), on peut relativiser cette notion comme suit : on dit
qu'un couple  $(Q,s)$ est un couplage ( distingu\'e ) {\bf partiel}
si $Q\subseteq P$ est filtrante d\'ecroissante, et $s$ est une
application $Q \rightarrow Q$ telle que sa restriction
constitue un couplage distingu\'e du sous-poset $Q$. On
note $Q=\dom(s)$. \\

 Soit ${\cal I} (P) $ l'ensemble des couplages partiels de $P$ ;
on a une relation d'ordre ${\leq}_{\cal I}$ naturelle sur 
 ${\cal I} (P) $, \`a savoir
$(Q_1,s_1) {\leq}_{\cal I} (Q_2,s_2)$ ssi
$Q_1 \subseteq Q_2$ et $s_2$ \'etend $s_1$. Nous appelons
{\bf couplages maximaux} les \'el\'ements qui sont
maximaux pour ${\leq}_{\cal I}$.
  On peut descendre d'un degr\'e encore dans la relativit\'e
en introduisant pour $Q \subseteq P$ filtrante d\'ecroissante
la notion de {\bf couplage Q-maximal de P} : les couplages
partiels $s$ de $P$ ayant cet \'epith\`ete sont ceux v\'erifiant 
la condition
$\forall s'$ \'etendant $s$, $\dom(s) \cap Q=\dom(s') \cap Q$.

 Si on suppose $Q$ finie, toute chaine finie
${\phi}_1,{\phi}_2, \ldots {\phi}_r$ avec chaque ${\phi}_{i+1}$
\'etendant ${\phi}_i$ et 
$\forall i, Dom({\phi}_i)\cap Q \subset Dom({\phi}_{i+1}) \cap Q$ 
(inclusion stricte) est n\'e\-cessairement de cardinal 
$\leq |Q|$ et si cette
chaine est de longueur maximale son dernier \'el\'ement
est un couplage $Q$-maximal ; donc : \newline

\begin{merkadenn}\label{unique_q_maximal_embedding}
\it Soit P un poset gradu\'e, Q une
partie filtrante d\'ecroissante finie de P. Alors tout couplage
partiel de $P$ se prolonge en un couplage $Q$-maximal. 
\end{merkadenn}

Commencons par donner un r\'esultat de ``passage du local au
global'' : \newline

\begin{teorem}\label{unique_p_maximal_embedding}
\it Soit P un poset gradu\'e tel que 
$\lbrace x \in P \ ; \ l(x)=k \rbrace$ soit fini pour tout $k$.
Si $A$ est une partie filtrante d\'ecroissante de $P$, alors
tout couplage partiel $\alpha$ d\'efini sur $A$ se prolonge
en un couplage maximal sur $P$. 
 Si de plus la fonction $coat$ est injective sur $P \setminus A$,
alors cette extension est unique.
\end{teorem}

 {\bf Preuve : } Existence.\newline
                
    Pour chaque $k$ notons 
$B_k=\lbrace x \in P \ ; \ l(x)\leq k \rbrace$.
Par la remarque 2.1, $\alpha$ admet une extension
${\phi}_0$ qui est $B_0$-maximale. En it\'erant
cette m\^eme remarque~\ref{unique_q_maximal_embedding}, 
on construit une suite
 $({\phi}_n)_{n \geq 0}$ de couplages partiels de $P$ tels
que 
\begin{displaymath}
\forall n \geq 1,
{\phi}_{n}  { \rm \ \acute{e}tend \ } {\phi}_{n-1}, \
{\phi}_{n} { \rm \ est \ } {B_n}-{ \rm maximale}.
\end{displaymath}
          
 Une fois cette suite $({\phi}_n)$ construite, d\'efinissons
$Q=\bigcup{Dom({\phi}_n)}_{n\geq 0}$, 
$\phi:Q \to Q$ par
\mbox{$\forall x \in Q, \phi(x)=\phi_n(x)$} si $x \in Dom( \phi_n)$. Alors
$\phi$ est bien d\'efinie et est un couplage
maximal \'etendant $\alpha$, comme cherch\'e.\newline

 Unicit\'e (dans le cas $coat$ injective sur $P \setminus A$).\newline

 Supposons par l'absurde que l'on ait deux couplages maximaux
distincts ${\mu}_1$ et ${\mu}_2$ qui \'etendent $\alpha$.
Prenons alors $w$ de longueur minimale tel que $\mu_1$ diff\`ere
de $\mu_2$ en $w$, i.e. (quitte \`a \'echanger $\mu_1$ et $\mu_2$)
\begin{displaymath}
\begin{array}{ll}
{ \rm Cas} \ 1 : &w \in Dom( \mu_1 ), \ w \in Dom( \mu_2 ), \ 
\mu_1(w) \neq \mu_2(w), {\rm ou \ bien} \\ 
{ \rm Cas} \ 2 : &w \in Dom( \mu_1 ), \ w \not\in Dom( \mu_2 ). 
\end{array}
\end{displaymath}

 Consid\'erons le cas 1. On a certainement $w\not\in A$ ;
et comme $\mu_1,\mu_2$ sont involutives, on doit
avoir $w \lhd \mu_1(w), \ w\lhd \mu_2(w), \
\mu_1(w) \not\in A, \ \mu_2(w) \not\in A$. Alors la condition
$(iii)$ de la d\'efinition d'un couplage maximal
donne $coat(\mu_1(w))=coat(\mu_2(w))$ puis
$\mu_1(w)=\mu_2(w)$ qui est exclu. \newline

 Passons au cas 2. De mani\`ere analogue au cas 1,
on voit que $w \not\in A, w \lhd \mu_1(w), \mu_1(w) \not\in A$.
On a m\^eme $\mu_1(w) \not\in Dom( \mu_2)$ (car $Dom( \mu_2)$ est
filtrante d\'ecroissante et $w\not\in Dom(\mu_2)$).
Consid\'erons $Q=Dom(\mu_2) \cup \lbrace w; \mu_1(w) \rbrace$
et $\phi \ : \ Q \to Q$ d\'efinie par
$\phi(x)=\mu_2(x)$ si $x\in Dom( \mu_2) $ et
$\phi(x)=\mu_1(x)$ si $x \in \lbrace w; \mu_1(w) \rbrace$. Alors 
$\phi$ est un couplage partiel qui \'etend strictement
$\mu_2$ ce qui contredit la maximalit\'e de ce dernier. \hfill {Q. E. D.}\\

 Nous allons maintenant quitter le monde des posets
g\'en\'eraux et abstraits pour nous restreindre jusqu'\`a
la fin de cet article au cas particulier ou $P$ provient
d'un syst\`eme de Coxeter $(W,S)$ de la mani\`ere suivante :
$P=W$, la relation d'ordre est l'ordre de Bruhat-Chevalley,
la longueur est la fonction longueur usuelle sur un groupe
de Coxeter.

 Dans ce cas particulier pr\'ecis, le th\'eor\`eme 
\ref{unique_p_maximal_embedding} prend
une forme plus fine. \'Etant donn\'e un
syst\`eme de Coxeter $(W,S)$, on appelle 
{\bf sous-groupe dih\'edral} de $W$ tout sous-groupe $<s,t>$ o\`u
$ \lbrace s;t \rbrace$ est une paire d'\'el\'ements de $S$.
Un \'el\'ement $ w \in W$ est dit 
{\bf \'el\'ement dih\'edral} si il appartient
\`a un sous-groupe dih\'edral, ou bien de mani\`ere
\'equivalente, si $w$ s'\'ecrit
\begin{displaymath}
w=\underbrace{sts \ldots}_{m \ termes}
\end{displaymath}
pour un certain entier $m$ et une certaine paire
$\lbrace s;t \rbrace \subseteq S$. Comme les \'el\'ements
dih\'edraux apparaissent \`a chaque instant dans ce travail, nous
introduisons tout de suite les notations suivantes :
 pour toute paire $\lbrace s;t \rbrace \subseteq S$ on note
\begin{displaymath}
\begin{array}{l}
[s,t,n \rangle=\underbrace{stst \ldots}_{n \ termes} \\
\langle n,t,s]=\underbrace{ \ldots tsts}_{n \ termes} \\
M_{st}=[s,t,m_{st}\rangle=\langle m_{st},t,s] \ {\rm si} \ m_{st}<\infty .
\end{array}
\end{displaymath}

Rappelons deux r\'esultats  d\'emontr\'es
ailleurs par Dyer et Waterhouse respectivement :

\begin{kinnig} \label{equal_coatoms_imply_equal_elements}
 \it Soit (W,S) un systeme de Coxeter.\newline
\begin{tabular}{l}
1) Pour $w\in W$, ($w$ est dih\'edral)$\Leftrightarrow$ 
($|coat(w)| \leq 2$). \\
2) Si $x$ et $y$ dans $W$ v\'erifient $coat(x)=coat(y)=A$ 
et $|A| \geq 3$, alors $x=y$.
\end{tabular}
\end{kinnig}

 {\bf Preuve :} Consulter [\ref{Dyer_1987}, proposition 7.25] pour la
premi\`ere assertion et [\ref{Waterhouse}, pro\-position7] pour la 
deuxi\`eme. \hfill {Q. E. D.}\\

 Cette proposition fait pressentir l'importance des \'el\'ements
dih\'edraux ; en fait, nous devons encore introduire la notion
d'{\bf \'el\'ement dih\'edral principal}\nolinebreak{:} 
si $(Q,\phi)$ est un couplage
partiel avec $Q \neq \emptyset$, alors $e \in Q$ et
par la r\`egle $(iii)$, $a=\phi(e)$ est un \'el\'ement de $S$.
Nous appelons {\bf sous-groupes dih\'edraux principaux}
les $P_s=<s,a>$ pour $s \in S \setminus \lbrace a \rbrace$ ; et 
les \'el\'ements dih\'edraux principaux sont les \'el\'ements
de 
\begin{displaymath}
P=\bigcup_{s\in S \setminus \lbrace a \rbrace}{P_s}.
\end{displaymath}

 Commen\c{c}ons par pr\'eciser l'action d'un couplage
maximal sur un sous-groupe dih\'edral : 

\begin{kinnig}\label{action_of_phi_on_dihedrals}
\it Soit $(W,S)$ un syst\`eme de
Coxeter, et $c=(Q,\phi)$ un couplage maximal sur $W$, $D$ un
sous-groupe dih\'edral de $W$. \newline
 (i) Si $D$ est principal, alors $\phi$ est d\'efini sur tout $D$
et $D$ est stable par $\phi$.\newline
 (ii) Si $D$ est non principal, alors
\begin{displaymath}
\forall w\in Q \cap D,\ { \rm on \ a \ } w \lhd \phi(w), \phi(w)\not\in D. 
\end{displaymath}
\end{kinnig}

 {\bf Preuve : } Soit $s \in S \setminus \lbrace a \rbrace$ et
$m=m(a,s)$ (coefficient entier ou infini de la matrice de Coxeter).
Rappelons que $P_s$ a un unique \'el\'ement de longueur $0$,
un ou pas d'\'el\'ement de longueur $m$ suivant
que $m$ est fini ou non,
et deux \'el\'ements en longueur $j$ pour $0<j<m$.
Pour $w\in {P_s}$ tel que $0<l(w)<m$, on notera 
$\bar{w}$ l'unique \'el\'ement de $P_s$ de m\^eme longueur
que, mais diff\'erent de, $w$. 

 Montrons $(i)$, c'est \`a dire 
\begin{displaymath}
\forall w\in P_s, \ w\in Q, \phi(w) \in P_s.
\end{displaymath}  
 
 On raisonne par r\'ecurrence sur $j=l(w)$. Si $j=0$,
on a $w=e$ donc $\phi(w)=a$ et le r\'esultat est clair.
Si $j=1$, on a $w=a$ (et alors $\phi(w)=e$) ou bien
$w=s$ (et alors $Z(\phi,w)=\lbrace a;s \rbrace$, 
donc $\phi(w) \in \lbrace as;sa \rbrace \subseteq P_s$).

 Soit maintenant $j \geq 2$. Supposons le 
r\'esultat vrai pour les longueurs $<j$.
Soit $w\in P_s$ tel que $l(w)=j$. Si $\phi(w) \lhd w$,
on a certainement $\phi(w) \in P_s$ car $P_s$ est 
filtrante d\'ecroissante. Sinon $w \lhd \phi(w)$
(ou bien $\phi(w)$ n'est pas d\'efini). En prenant un coatome
$v$ de $w$, on a $coat(w)=\lbrace v;\bar{v} \rbrace$ car $j\geq 2$,
donc $Z(\phi,w)=\lbrace w \rbrace \cup \phi(B)$ ou
l'on a pos\'e 
$B=\lbrace z \in \lbrace v;\bar{v}\rbrace \ ; \ z \lhd \phi(z)\rbrace$.

 Si $B=\lbrace v;\bar{v}\rbrace$ alors $\phi$ induirait
une bijection $\lbrace v;\bar{v}\rbrace \to \lbrace w;\bar{w}\rbrace$
ce qui est exclu car $w \lhd \phi(w)$. Ainsi
$B \neq \lbrace v;\bar{v}\rbrace$. 
  
 Supposons $B=\emptyset$. Alors en posant $u=\phi(v)$ on
a $\phi(v)=u \lhd v$, $\phi(\bar{v})=\bar{u} \lhd \bar{v}$.
De plus, comme $\phi$ est un couplage distingu\'e
on doit avoir $coat(\phi(u))=Z(\phi,u)$ c'est-\`a-dire
$coat(v)=\lbrace u \rbrace \cup \phi(B')$ ou l'on a pos\'e
$B'=\lbrace z \in coat(u) \ ; \  z \lhd \phi(z)\rbrace$.
Or $coat(v)=\lbrace u ; \bar{u}\rbrace$, donc
$\bar{u} \in \phi(B')$, donc $\phi(\bar{u}) \in B'$,
donc $l(\phi(\bar{u}))=l(u)-1$ qui contredit
$\bar{u} \lhd \phi(\bar{u})$. Ainsi $B\neq\emptyset$.
 
 Finalement $|B|=1$ et par exemple $B=\lbrace v \rbrace$.
Alors n\'ecessairement $\phi(v)=\bar{w}$ donc
$Z(\phi,w)=\lbrace w; \bar{w} \rbrace$. Par cons\'equent
$\phi(w)$ est d\'efini et est dans $P_s$, par la
proposition 2.3. Ceci ach\`eve la
preuve par r\'ecurrence et montre (i).\newline

 Montrons maintenant (ii).
Supposons par l'absurde que l'on ait $w \in D$ tel que
$\beta(w)\in D$ (cas 1) ou tel que $\beta(w) \lhd w$ (cas 2); on prend pour
$w$ un contre-exemple de longueur minimale. 
Supposons que $w$ est dans le cas 2. 
Alors $\beta(w)$ est
dans ce qu'on a appel\'e le cas 1, ce qui contredit la minimalit\'e de $w$.
  Ainsi $w$ est dans le cas 1, et pas dans le cas 2 : 
$w\lhd \beta(w), \beta(w)\in D$.
 Remarquons que $w \neq e$ car $\phi(e)=a$.
 Soit $u$ un coatome de $w$ ; alors par minimalit\'e de $w$,
$u \lhd \beta(u), \beta(u) \not\in D$,  donc $\beta(u)\in Z(\phi,w)$.
Mais alors 
$\beta(u)\lhd \beta(w), \ \beta(u)\not\in D, \beta(w)\in D$ ce
qui est absurde. \hfill {Q. E. D.}\newline

\begin{kinnig}\label{equality_on_pp_suffices}
\it Soit $(W,S)$ un syst\`eme de
Coxeter, $\beta$ et $\gamma$ deux couplages maximaux
sur $W$, et $w\in W$. Supposons que l'on ait :
\begin{displaymath}
\forall v \in [e,w] \cap P, \ \beta(v)=\gamma(v)
\end{displaymath}
 Alors 
\begin{displaymath}
\begin{array}{c}
\forall w \in W, \beta(w)=\gamma(w) \\
({  ou \ bien \ } \beta(w) 
{  \ et \ } \gamma(w) 
{ \ sont \ tous \ deux \ non \ d\acute{e}finis }).
\end{array}
\end{displaymath}
\end{kinnig}

 {\bf Preuve : } Supposons que la proposition soit fausse ;
prenons alors un contre-exemple $w$ de longueur minimale.
En raisonnant comme \`a la partie ``unicit\'e" du th\'eor\`eme 
\ref{unique_p_maximal_embedding}, on voit
que ce contrexemple $w$ v\'erifie n\'ecessairement 
$\beta(w)\neq \gamma(w), \ w\lhd \beta(w), \ w \lhd \gamma(w), \ 
coat( \beta(w))=coat( \gamma(w))$. Par la proposition 
\ref{equal_coatoms_imply_equal_elements}, ceci implique 
que $\beta(w)$ et $\gamma(w)$ sont deux \'el\'ements d'un m\^eme sous-groupe
dih\'edral $D$. Ce sous-groupe ne peut \^etre principal puisque
par hypoth\`ese $\beta$ et $\gamma$ coincident sur $P\cap [e,w]$. Mais
alors \ref{action_of_phi_on_dihedrals}.(ii) 
montre que $\beta(w) \not\in D$ ce qui 
est absurde. \hfill {Q. E. D.}\newline 

 En combinant les deux propositions pr\'ec\'edentes on obtient : \newline
 
\begin{teorem}\label{mw_is_a_cartesian_product}
\it Soit $(W,S)$ un syst\`eme de
Coxeter.\newline
 (i) Pour tout couplage maximal $\phi$ sur $W$, chaque sous-groupe
dih\'edral principal $P_s$ est stable par $\phi$, d'o\`u un
couplage induit ${\phi}_s$  sur $P_s$. \newline
 (ii) R\'eciproquement, pour toute famille ${({\alpha}_s)}_{s\neq a}$  
avec chaque ${\alpha}_s$ un couplage sur $P_s$ tel que 
${\alpha}_s(e)=a$, on a un unique couplage maximal $\phi$ qui
\'etend la r\'eunion des ${\alpha}_s$ : $\forall 
s, \ {{\phi}_{|P_s}}={\alpha}_{s}$.
\end{teorem}

 {\bf Preuve : } Le (i) n'est bien s\^ur qu'une r\'ep\'etition de
\ref{action_of_phi_on_dihedrals}.(i).

 Montrons maintenant le resultat d'extension unique. Tout
d'abord, $P=\bigcup{P_s}$ est
filtrant d\'ecroissant, d`o\`u un
couplage partiel $\alpha:P \to P$ ; ceci nous donne
d\'eja l'existence de $\phi$ par la partie ``existence'' du 
th\'eor\`eme \ref{unique_p_maximal_embedding}. L'unicit\'e d\'ecoule de
la proposition \ref{equality_on_pp_suffices}. \hfill {Q. E. D.}\newline

\bigskip

 \section{Formules de descente pour
les polyn\^omes {\it R}}\label{exit_r}

 \hed Comme annonc\'e dans l'introduction,
nous utilisons librement ici les 
propri\'et\'es \'el\'ementaires des
polyn\^omes $R$, nous
r\'eferant \`a 
[\ref{Humphreys}] pour toute explication
suppl\'e\-mentaire.\newline

 Pour $s$ dans $S$ on note $L_s= \lbrace w \in W ; l(sw)>l(w)   \rbrace $.
Les formules suivantes qui donnent lieu \`a une m\'ethode
de calcul des polyn\^omes $R_{u,v}$ sont bien \mbox{connues :}
 pour $x,y \in L_{s}, $ on a\newline

\hspace*{4cm}{$R_{sx,sy}=R_{x,y}$} \hfill (3.1.1)\\
\hspace*{3.5cm}{$R_{x,sy}=(q-1)R_{x,y}+qR_{sx,y}$}\hfill (3.1.2)\\

 L'id\'ee est d'essayer de montrer que ces formules restent valables,
{\it mutatis mutandis}, lorsque l'on remplace $s$ par
un couplage distingu\'e quelconque.\newline

\begin{termenadenn}\label{full_elements} 
\it  Soit (W,S) un syst\`eme
de Coxeter et $w \in W$. On dit que $w$ est {\bf plein par
rapport \`a} $J \subseteq S$ si 
$\forall s,t \in J, \ m_{st}< \infty, \  M_{st}\leq w$. 
On dit que $w$ est {\bf plein} s'il est $plein$ par rapport
\`a $S$. 
\end{termenadenn}

\begin{lemma}\label{subgroups_remain_full}
\it Soit $(W,S)$ un syst\`eme de Coxeter,
$J \subseteq S$, $w\in W$ plein
par rapport \`a $J$. Alors il existe $v \in <J>$ plein
par rapport \`a $J$ tel que $v \leq w$.
\end{lemma}

 {\bf Preuve : } On sait que $[e,w] \cap <J>$ a un plus
grand \'el\'ement $v$ (cf. par exemple [3, proposition 2.5]) ;
montrons que $v$ r\'epond \`a la question. On a certainement
$v \leq w$. De plus, si $s$ et $t$ sont
deux \'el\'ements distincts de $J$ et $\mu=M_{st}$ 
l'\'el\'ement dih\'edral maximal associ\'e, on
a $\mu \in [e,w] \cap <J>$ donc
$\mu \leq v$. Ceci \'etant vrai pour tous les
$s$ et $t$, $v$ est plein par rapport 
\`a $<J>$. \hfill{Q. E. D.}\newline

 Nous avons maintenant besoin des ensembles de descente
\`a gauche et \`a droite d'un \'el\'ement $w$ de $W$ : ce
sont respectivement 
$\lbrace s \in S \ ; \ sw \lhd w \rbrace$ et
$\lbrace s \in S \ ; \ ws \lhd w \rbrace$. On les note
$D_g(w)$ et $D_d(w)$ dans la suite de cet article.\newline

\begin{termenadenn}\label{regularity}
 \it Si $s \in S$ et $c=(Q, \phi)$ est un couplage,
on dit que $s$ est $c${\bf-r\'egulier} (\`a gauche)
ou que $c$ est $s$-r\'egulier (en toute rigueur la
notion de r\'egularit\'e concerne le couple $(s,c)$) si

\begin{displaymath}
\begin{array}{lll} 
(i) Q \ est \ stable \ par \ {(x \mapsto sx)}, \\
et\\
(ii) \forall x \in Q, \ \phi (sx)=s \phi (x).
\end{array}
\end{displaymath}
\end{termenadenn}
 Bien entendu, on a une d\'efinition analogue \`a droite.

\begin{termenadenn}\label{reducibility}
\it Soit (W,S) un syst\`eme
de Coxeter et $c=(Q,\phi)$ un couplage 
maximal de W ; soit $\omega$ une orbite dans $Q$
pour l'action de l'involution $\phi$. Alors
$\omega$ peut s'\'ecrire
$\omega=\lbrace m, M \rbrace$ 
avec $m \lhd M, \ \phi(m)=M$.
L'orbite $\omega$ est dite
{\bf pleine} si $M$ est plein.
On dit que $\omega$ est une
{\bf orbite r\'eductible \`a gauche} pour $c$
si il existe $s$ r\'egulier \`a gauche
dans l'ensemble de descente \`a gauche
de $m$. De m\^eme, $\omega$ est une
{\bf orbite r\'eductible \`a droite} pour $c$
si il existe $s$ r\'egulier \`a droite
dans l'ensemble de descente \`a droite
de $m$. L'orbite $\omega$ est dite {\bf orbite r\'eductible}
si elle est r\'eductible \`a gauche ou \`a droite.
Enfin, $c$ est un {\bf couplage r\'eductible} 
si $|S| \leq 2 $ ou si toute orbite pleine est
r\'eductible.
\end{termenadenn}

 On \'etend l'addition et la relation d'ordre
usuelle sur $\mathbb{N}=\lbrace 0;1;2; \ldots \rbrace$
\`a $\mathbb{N} \cup \lbrace \infty \rbrace$ par 
$x \leq \infty$ et $x+\infty=\infty$ pour
$x \in \mathbb{N} \cup \lbrace \infty \rbrace$. Le 
r\'esultat essentiel de cet section s'\'enonce
ainsi : \newline

\begin{kinnig}\label{descent_formulas}
\it Soit (W,S) un syst\`eme
de Coxeter ayant la propri\'et\'e suivante :
pour tout syst\`eme de Coxeter $(W',S)$ 
associ\'e \`a une matrice de Coxeter $M'$
v\'erifiant $\forall s,t \in S, \ m'_{st} \leq m_{st}$,
on a que tout couplage maximal de $W'$ est r\'eductible. 
Soit alors $c=(Q,\phi)$ un couplage distingu\'e de W. 
D\'efinissons l'ensemble $L_{\phi}= \lbrace w \in Q \ 
; \ w \lhd \phi (w)  \rbrace $. Alors, 
pour $(x,y) \in {L_{\phi}}^2, $ on a\\

\hspace*{4cm}{$R_{\phi(x),\phi(y)}=R_{x,y}$} 
\hfill (\ref{descent_formulas}.1)\\
\hspace*{3.5cm}{$R_{x,\phi(y)}=(q-1)R_{x,y}+qR_{\phi(x),y}$}
\hfill (\ref{descent_formulas}.2)\\
\end{kinnig}

  {\bf Preuve de la proposition.} Tout d'abord, remarquons
que quand $|S|\leq 2$, on a l'\'equivalence
$(u<v) \Leftrightarrow (l(u)< l(v))$, d'o\`u
on d\'eduit assez facilement que $R_{u,v}$ ne d\'epend que de
$l(v)-l(u)$ (on peut par exemple d\'efinir une suite de
polyn\^omes $(L_i(q))$ par $L_0=1, \ L_1=q,\  
\forall n \geq 2 \ L_n=(q-1)L_{n-1}+qL_{n-2} $ et montrer
$\forall u \leq v, R_{u,v}=L_{l(v)-l(u)}(q)$ par r\'ecurrence
sur la longueur de $v$ en utilisant \ref{descent_formulas}.1 
et \ref{descent_formulas}.2). On obtient
donc le r\'esultat tr\`es vite dans ce cas $|S| \leq 2$. \newline

 Ensuite, on raisonne par r\'ecurrence sur $y$ et
sur la taille du groupe de Coxeter~: formellement,
l'ensemble $S$ est fix\'e et on montre une propri\'et\'e
du couple $(M,y)$ en raisonnant par r\'ecurrence sur
la quantit\'e $q(M,y)=l(y)+||M||$, o\`u l'on pose
$||M||=\sum_{s,t \in S}{m_{st}}$ (\`a priori le
raisonnement par r\'ecurrence ne montre le r\'esultat
que pour les couples $(M,y)$ tels que
$q(M,y)$ soit fini~; mais il est facile
de voir que le raisonnement
de r\'eduction au cas ``$\phi(y)$ plein''
que nous allons exposer permet \'egalement de d\'eduire
le cas $q(M,y)=\infty$ du cas
$q(M,y)$ fini).\newline

 On va montrer que l'on peut
se ramener au cas $w=\phi(y)$ plein. En effet, supposons
$w$ non plein et consid\'erons
la matrice de Coxeter $M'$ d\'efinie par
${m'}_{st}=$la longueur du plus grand \'el\'ement
de $[e,w] \cap <s,t>$ pour $s,t \in S$ et
le syst\`eme de Coxeter $(S,W')$
associ\'e \`a la matrice $M'$. Par la proposition 3.5. de [3],
l'application (o\`u les $s_1 \ldots s_r$ sont
des mots r\'eduits)
\begin{displaymath}
\begin{array}{lll}
\psi \ : \ [e,w]& \to& W' \\
{ \lbrace s_1 \ldots s_r \rbrace}_{W} & \mapsto & 
{ \lbrace s_1 \ldots s_r \rbrace}_{W'}
\end{array}
\end{displaymath}
est bien d\'efinie, strictement croissante pour les ordres de Bruhat
et v\'erifie
\begin{displaymath}
\forall x \in [e,w], \ \forall s \in S,
(xs \in [e,w]) \Rightarrow (\psi({\lbrace xs \rbrace}_{W})=
{ \lbrace \psi(x)s \rbrace}_{W'})
\end{displaymath}
 De ceci on d\'eduit facilement que $\psi$ r\'ealise un
isomorphisme de posets gradu\'es de $[e,w]$ sur
$[e,\psi(w)]$ et que
\begin{displaymath}
\forall u,v \in [e,w], \
R^{W}_{u,v}=R^{W'}_{\psi(u),\psi(v)}
\end{displaymath}

 (raisonner par r\'ecurrence sur la longueur de $v$ en
utilisant les formules \ref{descent_formulas}.1. et 
\ref{descent_formulas}.2). 
Par cons\'equent, le probl\`eme sur
$[e,w] \subseteq W$ se transporte compl\`etement
sur $[e,\psi(w)] \subseteq W'$ dans lequel 
effectivement $\psi(w)$ est plein. Comme $w$ est
non plein on a $||M'||<||M||$ donc 
$q(M',\psi(y))<q(M,y)$ d'o\`u le r\'esultat par
r\'ecurrence lorsque $w$ est non plein. \newline 

 \`A partir de maintenant, on reste dans un
groupe de Coxeter fix\'e dans lequel
$\phi(y)$ est plein ; en particulier
$||M||< \infty$. La r\'ecurrence sur $q(M,y)$
se r\'eduit alors simplement \`a une r\'ecurrence
sur la longueur de $y$.\newline

 Cas $l(y)=0$ :\newline
 \hspace*{.5cm}Dans ce cas $y=e$ et toutes les $R_{u,v}$ 
consider\'es sont nuls sauf  si $x=e$ ou $\phi(e)$ ; dans
chacun de ces cas, on v\'erifie directement les
formules \ref{descent_formulas}.1 et \ref{descent_formulas}.2.\newline

 Cas $l(y)>0$ avec le r\'esultat vrai pour 
les $y'$ de longueur $<l(y)$ :\newline
 \hspace*{.5cm} Gr\^ace aux \'egalit\'es bien
connues $R_{u,u}=1$ et $R_{u,v}=q-1$ si $u \lhd v$,
on peut supposer $x<y$.

Supposons par exemple qu'il existe $g$ r\'egulier
\`a gauche dans l'ensemble de
descente \`a gauche  de $y$, le cas \`a droite
\'etant tout-\`a-fait sym\'etrique.\\

 Soit $v=gy$. Si $m=l(v)$, on a donc $l(y)=m+1$,
$l(\phi(y))=m+2$. Si $p=l( \phi(v))$, on a
d'une part $p-l(v) \in \lbrace -1,1 \rbrace$ et d'autre part
$p-l(g \phi(v)) \in \lbrace-1,1 \rbrace$ donc $p=m+1$,
et finalement $v \lhd gv \lhd \phi(y)$, 
$v \lhd \phi(v) \lhd \phi(y)$.\newline

 Supposons d'abord $gx \lhd x$. Alors, en posant
$w=gx$, le raisonnement qui vient d'\^etre fait
(avec $x$ au lieu de $y$) donne
$w \lhd gw \lhd \phi(x)$, $w \lhd \phi(w) \lhd \phi(x)$.
Dans ce cas,
\begin{trivlist}
\item[]
\hspace*{1.5cm}$R_{\phi (x), \phi (y) }=R_{g\phi(w), g\phi (v)}
=R_{\phi(w), \phi(v)}$

\item[]
\hspace*{3cm}$=R_{w,v}$ \hfill (hypoth\`ese de r\'ecurrence)
\item[]                        
\hspace*{3cm}$=R_{gw,gv}=R_{x,y}.$
\end{trivlist}
et
\begin{trivlist}
\item[]
\hspace*{2cm}$R_{x, \phi (y) }=R_{gw, g\phi (v)}
=R_{w, \phi(v)}$

\item[]
\hspace*{3cm}$=(q-1)R_{w,v}+
qR_{\phi(w),v}$ \hfill (hypoth\`ese de r\'ecurrence)
\item[]                        
\hspace*{3cm}$=(q-1)R_{gw,gv}+qR_{g\phi(w),gv}=
(q-1)R_{x,y}+qR_{\phi(x),y}.$
\end{trivlist}

 Traitons maintenant le cas $x \lhd gx$. Alors, comme
on a $coat(\phi(x))=\lbrace x \rbrace \cup 
\lbrace \phi(z) \ ; \ z \lhd x, z \lhd \phi(z) \rbrace$, on voit
que $\phi(gx)$ ne peut \^etre un coatome de $\phi(x)$ que
si $\phi(gx)=x$, i.e. si $\phi(x)=gx$ ; sinon
on a $\phi(x) \lhd \phi(gx)$. On a alors deux
\mbox{sous-cas :}
\begin{enumerate}
\item[] $x \lhd gx$, $\phi(x)=gx$,
\item[] $x \lhd gx$, $\phi(x) \lhd \phi(gx)$.
\end{enumerate}

 Dans le premier sous-cas, on a
\begin{trivlist}
\item[] \hspace*{1.5cm}$R_{\phi(x),\phi(y)}=R_{gx,g\phi(v)}=R_{x,\phi(v)}$
\item[] \hspace*{3cm}
$=(q-1)R_{x,v}+qR_{\phi(x),v}$ \hfill (hypoth\`ese de r\'ecurrence)
\item[] \hspace*{3cm}
$=(q-1)R_{x,v}+qR_{gx,v}=R_{x,gv}=R_{x,y}$.
\end{trivlist}
et
\begin{trivlist}
\item[] $R_{x,\phi(y)}=R_{x,g\phi(v)}=
(q-1)R_{x,\phi(v)}+qR_{gx,\phi(v)}=
(q-1)R_{x,\phi(v)}+qR_{\phi(x),\phi(v)}$
\item[] \hspace*{1.5cm}
$=(q-1)((q-1)R_{x,v}+qR_{\phi(x),v})
+qR_{x,v}$ \hfill (hypoth\`ese de r\'ecurrence)
\item[] \hspace*{1.5cm}
$=(q-1)((q-1)R_{x,v}+qR_{gx,v})
+qR_{gx,gv}$
\item[]\hspace*{1.5cm}
$=(q-1)R_{x,gv}+qR_{gx,gv}=
(q-1)R_{x,y}+qR_{\phi(x),y}$.
\end{trivlist}

 Et pour finir, dans le deuxi\`eme sous-cas on a

\begin{trivlist}
\item[]
\hspace*{3cm}$R_{\phi (x), \phi (y) }=R_{\phi(x), \phi (gv)}
=R_{\phi(x), g \phi(v)}$
\item[]
\hspace*{4.5cm}                       $=(q-1)R_{\phi(x),\phi (v)}
                        +qR_{g \phi(x),\phi (v)}$

\item[]
\hspace*{4.5cm}                        $=(q-1)R_{x,v}
                        +qR_{gx,v}$ \hfill (hypoth\`ese de r\'ecurrence)
\item[]                        
\hspace*{4.5cm}       $=R_{x,gv}=R_{x,y}.$
\end{trivlist}
et
\begin{trivlist}
\item[]
$R_{x, \phi (y) }=R_{x, \phi (gv)}
=R_{x, g \phi(v)}$
\item[]
\hspace*{1cm}                       $=(q-1)R_{x,\phi (v)}
                        +qR_{gx,\phi (v)}$

\item[]
\hspace*{1cm}    $=(q-1)\lbrace (q-1)R_{x,v}+qR_{\phi(x),v} \rbrace
    +q \lbrace (q-1)R_{gx,v}+qR_{\phi(gx),v} \rbrace $
\item[]
\hspace*{1cm}    $=(q-1)\lbrace (q-1)R_{x,v}+qR_{gx,v} \rbrace
    +(q-1) \lbrace (q-1)R_{\phi(x),v}+qR_{g\phi(x),v} \rbrace $
\item[]
\hspace*{1cm}                       $=(q-1)R_{x,gv}
                        +qR_{\phi(x),v}$
\item[]
\hspace*{1cm}                       $=(q-1)R_{x,y}
                        +qR_{\phi(x),y}$. \hfill Q. E. D.

\end{trivlist}

 Notre but va consister maintenant a montrer
que tous les couplages sont r\'eductibles.

\begin{termenadenn}\label{full_matchings}
\it  Soit (W,S) un syst\`eme
de Coxeter et $c=(Q,\phi)$ un couplage sur W. 
On dit que $c$ est {\bf plein} si $Q$ contient un \'el\'ement
plein. 
\end{termenadenn}

 Remarquons qu'un couplage non plein est trivialement
r\'eductible. Ceci \'eliminera un bon nombre de cas
dans ce qui va suivre. \newline

 \section{Crit\`eres de r\'egularit\'e}\label{big_enough}

 \hed Ces crit\`eres de r\'egularit\'e permettent
 non seulement de montrer
qu'un couplage est d\'efini en certain points, mais
en plus donnent une formule explicite pour la valeur du
couplage en ces points. Dans ce qui suit,
``$\phi(x)=y$ '' sous-entend ``$\phi$ d\'efinie en 
$x$ et $\phi(x)=y$.''\newline

\begin{kinnig}\label{commutativity_on_ps_suffices}   
\it Soit (W,S) un syst\`eme de
Coxeter, $c=(Q,\phi)$ un couplage maximal sur W,
$a=\phi(e)$. Soit $w\in Q$
et $s\in S$.
(on rappelle que pour $s\neq a$ on pose $P_s=<s,a>$
et $P=\bigcup_{s \neq a}{P_s}$)
\begin{displaymath}
\begin{array}{ll}
Si \ s\neq a, \ et \ \forall v\leq w, \  
(v \in P_s) \Rightarrow (\phi(sv)=s\phi(v)), \ alors \
 \phi(sw)=s\phi(w).\\
Si \ s=a, \ et \ \forall v\leq w, \  
(v \in P) \Rightarrow (\phi(sv)=s\phi(v)), \ alors \ \
 \phi(sw)=s\phi(w).
\end{array}
\end{displaymath}
Bien entendu, on a une variante en faisant agir $s$
\`a droite.
\end{kinnig}

 {\bf Preuve : } Nous nous contenterons de d\'emontrer la
premi\`ere assertion, la deuxi\`eme \'etant tout-\`a-fait
similaire.\newline

  Supposons $u=\phi(w) \lhd w$. Alors par
hypoth\`ese de r\'ecurrence on a $\phi(su)=s\phi(u)=sw$,
donc l'assertion est  vraie dans ce cas.

 On raisonne par r\'ecurrence sur la longueur de $w$. Comme
il arrive souvent, le cas $l(w)=0$ est trivial, de m\^eme
d'ailleurs que le cas $w\in P_s$. On prend donc $w\not\in P_s$.
Consid\'erons la formule :
\begin{displaymath}
F(w)\ : \ \phi(sw)=s\phi(w).
\end{displaymath}
 
 Comme $\phi$ et $x \mapsto sx$ sont des involutions,
on a $F(w) \Leftrightarrow F(sw)$ et $F(w) \Leftrightarrow F(\phi(w))$.
Par cons\'equent, on peut supposer $w\lhd sw, w \lhd \phi(w) $.
Supposons $y=\phi(sw)\lhd sw$. Si $l(sy)>l(y)$, comme
$(y\lhd sw, \ w\lhd sw)$ cela
implique $y=w$, donc $\phi(w)=sw$ et l'assertion est claire.
Sinon on a $y=sz$ avec $l(y)=l(z)+1, z \lhd w.$ Par
hypoth\`ese de r\'ecurrence on a $F(z)$, donc $F(sz)=F(y)$,
donc $F(\phi(y))=F(sw)$, donc $F(w)$.
 Supposons $v=s\phi(w) \lhd \phi(w)$.
Alors on a $v=w$ (auquel cas on retrouve le cas
$\phi(w)=sw$ ) ou bien il existe $x\lhd w$ avec
$v=\phi(x),x \lhd \phi(x)$.  Alors par
hypoth\`ese de r\'ecurrence on a $F(x)$,
donc $F(\phi(x))=F(v)$, donc $F(w)$.\newline

 On peut donc supposer $w \not\in P_s, w \lhd sw,
w \lhd \phi(w), 
sw \lhd \phi(sw)$ (ou $\phi$ non d\'efini en $sw$) , 
$\phi(w) \lhd s\phi(w) $. \newline

 Si $s\phi(w)$ est dih\'edral, il existe un sous-groupe
dih\'edral $D$ tel que $s\phi(w) \in D$. Alors
$w \in D$ et $\phi(w) \in D$.
La proposition \ref{action_of_phi_on_dihedrals}
montre alors que $D$ est principal:
il existe $t\in S\setminus \lbrace a \rbrace$ 
tel que $D=P_t$.
Alors $s \in P_t$ 
donc $s=t$ et enfin $w\in P_s$ qui est exclu.
 
 Donc $s\phi(w)$ est non dih\'edral.

 Par la proposition \ref{equal_coatoms_imply_equal_elements}, 
tout $x$ v\'erifiant
$coat(x)=coat(s\phi(w))$ est en fait confondu avec
$s\phi(w)$. Ceci va nous permettre de montrer
une \'egalit\'e d'\'el\'ements via une \'egalit\'e
d'ensembles de coatomes.\newline

 On est en droit d'\'ecrire (en utilisant l'hypoth\`ese
de r\'ecurrence \`a la quatri\`eme ligne , et sous r\'eserve d'existence
pour le premier ensemble de coatomes)
\begin{displaymath}
\begin{array}{lll}
coat(\phi(sw))&=&
\lbrace sw \rbrace \cup 
\lbrace \phi(z) \ ; \ z \lhd sw, z \lhd \phi(z) \rbrace \\
&=&
\lbrace sw \rbrace \cup 
\lbrace \phi(z) \ ; \ (z=w \ { \rm ou} \ 
z=su, u \lhd w, u\lhd su), z \lhd \phi(z) \rbrace \\
&=&
\lbrace sw;\phi(w) \rbrace \cup 
\lbrace \phi(su) \ ; \ 
 u \lhd w, u\lhd su, su \lhd \phi(su) \rbrace \\
&=&
\lbrace sw;\phi(w) \rbrace \cup 
\lbrace s\phi(u) \ ; \ 
 u \lhd w, u\lhd su, su \lhd s\phi(u) \rbrace \\
{\rm et} \ coat(s\phi(w))&=&
\lbrace \phi(w) \rbrace \cup 
\lbrace sz \ ; \ z \lhd \phi(w), z \lhd sz \rbrace \\
&=&
\lbrace \phi(w) \rbrace \cup 
\lbrace sz \ ; \ (z=w \ { \rm ou} \ 
z=\phi(u), u \lhd w, u\lhd \phi(u)), z \lhd sz \rbrace \\
&=&
\lbrace sw;\phi(w) \rbrace \cup 
\lbrace s\phi(u) \ ; \ 
 u \lhd w, u\lhd \phi(u), \phi(u) \lhd s\phi(u) \rbrace \\

\end{array}
\end{displaymath}

 Pour conclure, il ne nous reste plus qu'\`a montrer que\newline
\centerline{$A= \lbrace u \ ; \ u \lhd w, u\lhd su, su \lhd s\phi(u) \rbrace$
et}\newline 
\centerline{$B= \lbrace u \ ; \ u \lhd w, u\lhd \phi(u), 
\phi(u) \lhd s\phi(u) \rbrace$}\newline
sont confondus. Par sym\'etrie des r\^oles de $\phi$ et $x \mapsto sx$, il
suffit de montrer $A \subseteq B$.

 Soit donc $u \in A$. Posons  $m=l(u)$ et $x=\phi(u)$.
On a donc $l(su)=m+1, l(s\phi(u))=m+2$.
Comme on a \`a la fois $l(x)-l(u) \in \lbrace -1;+1 \rbrace$
et $l(s\phi(u))-l(x) \in \lbrace -1;+1\rbrace$, ceci impose
$l(x)=m+1$ donc $u\in B$. \hfill {Q. E. D.}\newline

  Si l'on se souvient de la proposition 
\ref{mw_is_a_cartesian_product}, la proposition 
\ref{commutativity_on_ps_suffices}
plus haut montre que

\begin{heuliadenn}\label{regularity_on_ps_suffices} 
\it Soit (W,S) un syst\`eme de
Coxeter, $c=(Q,\phi)$ un couplage maximal sur W,
$a=\phi(e)$. Soit $w\in Q$
et $s\in S$.

\begin{displaymath}
\begin{array}{ll}
Si \ s\neq a, ( \phi \  est \ s-r\acute{e}gulier \ \grave{a} \ gauche ) 
\Leftrightarrow 
( {\phi}_{|P_s} \  est \ s-r\acute{e}gulier \ \grave{a} \ gauche )\\
Si \ s=a, ( \phi \  est \ s-r\acute{e}gulier \ \grave{a} \ gauche ) 
\Leftrightarrow 
( {\phi}_{|P} \  est \ s-r\acute{e}gulier \ \grave{a} \ gauche ).
\end{array}
\end{displaymath}
Avec bien \'evidemment une variante \`a droite.
\end{heuliadenn}

\begin{teorem}\label{explicit_formula_for_phi}
\it Soit (W,S) un syst\`eme de Coxeter, avec W muni
d'une structure de poset gradu\'e par l'ordre de Bruhat et
la fonction longueur usuelle. Soit $c=(Q, \phi)$ un couplage
maximal sur $W$, $a=\phi(e)$, $\lambda$ l'application
$W \to W, x \mapsto xa$,
$X$ et $Y$ deux parties de $S$
telles que :
\begin{displaymath}
\begin{array}{l}
\forall x \in X \setminus \lbrace a \rbrace, \ 
{\phi}_{|P_x}={\lambda}_{|P_x} \\
\forall y \in Y \setminus \lbrace a \rbrace, \ 
{\phi}_{|P_y} \ est \ a-r\acute{e}gulier \ \grave{a} \  gauche. \\
\end{array}
\end{displaymath}
Alors $<X>(<Y> \cap Q) \subseteq Q$ et 
\begin{displaymath}
\forall x \in <X>,
\forall y \in <Y> \cap Q, \ 
\phi(xy)=x\phi(y).
\end{displaymath}

 Bien entendu, on a un r\'esultat analogue 
en \'echangeant gauche et droite.
\end{teorem}

 {\bf Preuve du th\'eor\`eme :} On raisonne par r\'ecurrence sur
la longueur de $x$ pour montrer que 
\begin{displaymath}
\forall x \in <X>, \ H(x) \ :  \ 
\forall y \in <Y> \cap Q, \
 \phi (xy)=x \phi (y).
\end{displaymath}

 Clairement $H(e)$ est vraie. Supposons $H(x^{'})$
vraie pour les $x^{'}$ de longueur $<l(x)$. 
Prenons $x_1 \in X$ tel que l'on puisse \'ecrire
$x=x_1\xi$ avec $\xi \lhd x$. Notons 
$w=\xi y$. Alors par hypoth\`ese
de r\'ecurrence $\phi(w)=\xi \phi(y)$. On cherche
\`a montrer $\phi(xy)=x\phi(y)$, c'est-\`a-dire
$\phi(x_1w)=x_1\phi(w)$. Pour cela, gr\^ace au
th\'eor\`eme \ref{commutativity_on_ps_suffices}, il suffit de montrer :
\begin{displaymath}
\hspace*{3.5cm}
\forall v \in [e,w] \cap P_{x_1}, \ \phi(x_1 v)=x_1 
\phi(v). \hspace*{3cm}(*)
\end{displaymath}
(o\`u l'on pose par commodit\'e $P_a=P$). Prenons donc un
$v$ v\'erifiant $v \in [e,w] \cap P_{x_1}$.\newline

 Supposons $x_1\neq a$. Alors par hypoth\`ese
$\phi$ est confondue avec $\lambda$ sur $P_{x_1}$,
donc $\phi(v)=va, \phi(x_1v)=x_1va$ et (*) est clair.\newline

 Supposons $x_1=a$. Alors il existe
$t\in (X \cup Y) \setminus \lbrace a \rbrace$ tel que
$v \in P_t$. Par hypoth\`ese $\phi$ est $a$-r\'eguli\`ere \`a
gauche sur $P_t$ donc $(*)$ est vraie. \hfill {Q. E. D.} \newline
 
 Donnons pour finir un crit\`ere pratique
de r\'egularit\'e : \newline

\begin{merkadenn}\label{regularity_criterion}
\it Soit $\phi$ un couplage
d\'efini sur un groupe de Coxeter dih\'edral $<x,y>$. Pour
$j \leq m_{xy}$ notons 
\begin{displaymath}
Y_j=[y,x,j \rangle
\end{displaymath} 
Alors on a \'equivalence entre
\begin{displaymath}
\begin{array}{ll}
(1) & \phi \ n'est \ pas \ x-r\acute{e}gulier \ \grave{a} \ gauche \\
(2) & \exists j \leq m_{xy}-3, \ 
\phi(Y_j)=Y_{j+1}, \ \phi(xY_j)\neq xY_{j+1}(donc \ \phi(xY_j)=Y_{j+2}).
\end{array}
\end{displaymath}
\end{merkadenn}
 
 {\bf Preuve :} Notons $Z=\lbrace z \in <x,y> \ ; \ 
\phi(xz) \neq x\phi(z) \rbrace$. Comme $\phi$ et
$t \mapsto xt$ sont des involutions, $Z$ est stable
par ces deux applications. Par cons\'equent, tout
\'el\'ement minimal $z_0$ de $Z$, s'il existe, v\'erifie
$z_0 \lhd xz_0$ et $z_0 \lhd \phi (z_0)$ 
(ce qui implique $(2)$ avec $z_0=Y_j$),
d'o\`u le r\'esultat. \hfill {Q. E. D.}\newline

 \section{R\'eduction du domaine dans le cas crois\'e}\label{not_too_big}

 \hed \'Etant donn\'e un couplage maximal $(Q,\phi)$ sur
un groupe de Coxeter $W$, il est facile de voir
que pour chaque $s\in S$ on a
$\phi(s) \in \lbrace as;sa \rbrace$ (o\`u $a=\phi(e)$). 
Lorsque la restriction de $\phi$ aux g\'en\'erateurs
n'est ni une multiplication \`a gauche ni une multiplication
\`a droite, o\`u de mani\`ere \'equivalente
si il existe $u,v \in S$ avec 
$m_{au}>2, \ m_{av}>2, \phi(u)=ua, \phi(v)=av$, on dit que
$\phi$ est crois\'e.
 En fait, nous allons montrer ici un r\'esultat vrai en toute
g\'en\'eralit\'e (le th\'eor\`eme \ref{first_result_on_smallness_of_q})
mais qui ne nous sera utile que dans le cas 
{\bf crois\'e} (proposition \ref{crossed_case}).\newline
 
 Soit $(W,S)$ un syst\`eme de Coxeter quelconque, et
$(Q, \phi)$ un couplage maximal associ\'e \`a $W$,
$a= \phi (e)$. Soit
$G$ et $D$ les parties de $S$
d\'efinies par
\begin{displaymath}
\begin{array}{ll}
G= \lbrace g \in S \ ; \ \phi (g)=ga \rbrace\\
D= \lbrace d \in S \ ; \ \phi (d)=ad \rbrace\\
\end{array}
\end{displaymath}

et $<G>$ et $<D>$ les
sous-groupes paraboliques associ\'es. Nous allons
d\'emontrer l'inclusion suivante :\\

\begin{teorem}\label{first_result_on_smallness_of_q}
$Q \subseteq <G><D>$.
\end{teorem}

 {\bf Preuve du th\'eor\`eme. } Supposons par l'absurde 
qu'il existe $w$
dans \mbox{$Q \setminus <G><D>$}. On peut
prendre $w$ minimal ; alors 
\mbox{$\forall v<w, v \in <G><D>$}. Remarquons
d'abord que l'ensemble de descente \`a
gauche $D_g(w)$ de $w$ ne contient que des \'el'ements
qui ne sont pas dans $G$ (sinon on peut
\'ecrire $w=gv$ avec $g \in G$, $v<w$ et
alors $v$ est dans $<G><D>$ donc $w$ aussi:
impossible), et comme  $S=G \cup D$, ces
\'el\'ements sont dans $D \setminus G$. De
m\^eme, les \'elements de l'emsemble de
descente \`a droite de $w$ sont tous dans
$G \setminus D$ :
$D_d(w) \subseteq G \setminus D$.

Soit
$w_1 \ldots w_m$ une \'ecriture r\'eduite
de $w$. On a donc $w_1 \in D \setminus G$,
$w_m \in G \setminus D$. Par la remarque
pr\'ec\'edente, l'\'element $x$ de $W$ represent\'e
par $w_1 \ldots w_{m-1}$ v\'erifie
$D_g(x) \cap G= \emptyset$. Comme $x \in <G><D>$,
ceci impose $x \in <D>$. 
Ainsi (en utilisant la notion de support dans
un groupe de Coxeter), on a
$\forall i \leq m-1, w_i \in D$. De m\^eme,
$\forall i \geq 2, w_i \in G$. Ainsi, en
renommant les $w_i$,\\
\centerline{$w=d b_1 \ldots b_r g, \ {\rm avec}$}
\begin{displaymath}
\left(
\begin{array}{lll}
d \in D \setminus G,\\
\forall i \ b_i \in D \cap G,\\
g \in G \setminus D
\end{array}
\right) (*)
\end{displaymath}

 De plus, comme d'une part
$D_g(w) \subseteq \lbrace d;b_1; \ldots ; b_r;g\rbrace$
et d'autre part $D_g(w) \subseteq D \setminus G$,
on voit que $D_g(w)=\lbrace d \rbrace$, et de m\^eme
$D_d(w)=\lbrace g \rbrace$. Ainsi, toute \'ecriture
r\'eduite de $w$ comporte les caract\`eres $g$ et $d$ une
et une seule fois.\\

 Nous utilisons le r\'esultat (facile) suivant (cf [\ref{du_Cloux_2000}], corollaire
de la proposition 2.6 ):

\begin{merkadenn}\label{from_europe}
\it Si $q \in Q$ et a n'est pas dans
le support de q, alors $q\lhd \phi(q)$ et si $\mu$ est un
mot r\'eduit repr\'esentant $\phi (q)$ dans $W$, on
obtient une \'ecriture r\'eduite de $q$ en effa\c{c}ant
le caract\`ere $a$ de $\mu$.
\end{merkadenn}  

 Dans ce qui suit, on utilise assez souvent l'ensemble
$\cal I$ des \'el\'ements de $W$ qui ont une unique
\'ecriture r\'eduite; notamment, comme $d \in D \setminus G$,
$g \in G \setminus D$, les \'el\'ements $ag$, $da$, $gad$
et $dag$ sont dans $\cal I$. En fait, la seule propri\'et\'e
de $\cal I$ qui nous int\'eresse est la suivante :\\ 
\begin{center}
Si $f \in \cal I$ et $m_f$ est l'unique mot r\'eduit
repr\'esentant $f$, \\
alors pour tout $h \in W$ tel que $f \lhd h$ et tout 
mot r\'eduit $m_h$ \\
repr\'esentant
$h$, le mot $m_f$ est une sous-expression de $m_h$.
\end{center}
Remarquons maintenant que 

\begin{displaymath}
\begin{array}{ll}
(1) \ \ {\rm Si} \ dg \neq gd, \ dg \not\in Q  \\
(2) \ \ {\rm Dans \ tous \ les \ cas}, \  dag\not\in Q.
\end{array}
\end{displaymath}

 Pour montrer (1) et (2), on raisonne dans les deux cas
par l'absurde : si $dg \in Q$, la remarque 4.2. ci-dessus montre
que $dg \lhd \phi (dg)$, donc 
$coat(\phi(dg))= \lbrace dg ; ad ; ga \rbrace $, or aucun
\'el\'ement de $W$ n'a cet ensemble de coatomes
(si $coat(w)=\lbrace dg; ad; ga \rbrace$, comme
$dg \in \cal I$ et $dg \lhd w$ on a
$w=xdg, dxg$ ou $dgx$, avec $x$ un caract\`ere 
de $S$. Comme $a \geq w$, $x=a$, mais alors
$ad$ et $ga$ ne peuvent pas \^etre des coatomes de $w$
les deux \`a la fois), d'o\`u (1). Montrons maintenant (2), et
supposons donc $dag \in Q$. Par le (1) et le fait que
$Q$ est filtrant \`a gauche, on a $dg=gd$. On a alors
$\phi(da)\in \lbrace ada,dad \rbrace,
\phi(ag) \in \lbrace aga,gag \rbrace, \phi(dg)=gad$, donc
$\forall x \in coat(dag), \phi(x) \neq dag$. Ainsi
$dag \lhd \phi(dag)$ ; notons $w=\phi (dag)$. Alors, comme 
$dg \leq dag$ et $\phi(dg)=gad$, on
a $gad \lhd \phi(dag)$; comme $gad$ est
dans $\cal I$, si $m$ est un
mot r\'eduit repr\'esentant $w$, on a un caract\`ere
$x$ tel que
\begin{displaymath}
m \in \lbrace xgad,gxad,gaxd,gadx \rbrace
\end{displaymath}

 Comme $ag \lhd w$, et $ag \in \cal I$,
la seule possiblit\'e pour $m$
\`a la premi\`ere ligne est $m=xgad, x=a$.
De m\^eme $da \lhd w$ donne $m=gadx, x=a$.
Alors les mots $agad$ et $gada$ sont confondus,
ce qui est absurde, d'o\`u (2).

 Reprenant notre raisonnement de
d\'epart, ces faits (1) et (2) donnent $dg=gd$, et
$\forall i, b_i \neq a$
(sinon $dag \leq w$, ce qui est impossible car $Q$ est 
filtrant \`a gauche). 
Par
cons\'equent $a \not\leq w$, donc par la remarque
\ref{from_europe} ci-dessus $w \lhd \phi (w)$, et 
si $\mu$ est une \'ecriture r\'eduite de $\phi (w)$,  
$\mu^{\sharp}$ le mot obtenu en supprimant l'unique
occurrence de $a$ dans $\mu$, on a $w=\mu^{\sharp}$
dans $W$. Donc $\mu$
contient les caract\`eres $a$, $d$ et $g$ une et
une seule fois. Comme $d \leq w$ et $w \lhd \phi(w)$, 
on a $ad=\phi (d) \leq \phi(w)$ et de m\^eme
$ga \leq \phi (w)$, donc le seul ordre d'apparition
possible dans $\mu$ est $g,a,d$. On voit alors
que l'ordre d'apparition dans ${\mu}^{\sharp}$ est $g,d$ ce
qui est absurde. \hfill {Q. E. D.} \newline

\bigskip
\bigskip

\section{\'Etude partielle du cas des groupes de rang 3.} 
\label{rank_three}

 \hed Avant d'entrer dans le vif du sujet nous indiquons
quelques outils qui seront utilis\'es implicitement
sans plus d'explications par la suite. Les faits suivants
sont bien connus, pour $w$ \'el\'ement d'un 
groupe de Coxeter :
\begin{trivlist}
\item[(1)] \'Etant donn\'e deux \'ecritures r\'eduites
de $w$, on peut passer de l'une \`a l'autre
en utilisant uniquement des relations de tresses.
\item[(2)] \'Etant donn\'e une \'ecriture quelconque de
$w$, on peut aboutir \`a une \'ecriture r\'eduite en
utilisant uniquement les relations de tresses et
les relations $s^2=e$ pour $s\in S$.
\end{trivlist}
 
 Le (1) servira implicitement \`a justifier chaque assertion
du type ``tel \'el\'ement $w$ a une unique \'ecriture
r\'eduite''; plus pr\'ecisement, en utilisant la notation
$\cal I$ d\'efinie au \S4, on a pour tout mot r\'eduit $m$ repr\'esentant
un \'el\'ement $w\in W$, $w\in \cal I$ ssi aucune relation de tresse
n'est utilisable sur $m$ i.e. ssi $m$ ne contient pas de sous-mot
dih\'edral correspondant \`a un \'el\'ement dih\'edral maximal. 
 Tandis que (2) sera implicitement utilis\'e
chaque fois qu'on aura besoin de
savoir qu'un mot est r\'eduit. Remarquons que les mots
rencontr\'es ne seront jamais bien complexes (ils ne diff\`ereront
d'un mot dih\'edral que par au plus un caract\`ere), ce qui
justifie que nous ne nous y attardions pas.\newline

 Dans toute cette section 6, on consid\`ere un syst\`eme de
Coxeter $(W,S)$ de rang 3 : $S=\lbrace a;b;b'\rbrace$ et
$(Q,\phi)$ un couplage maximal sur $W$ avec $\phi(e)=a$.\newline

 \subsection{G\'en\'eralit\'es en rang 3.}\label{r3_general}

\begin{kinnig_lechel}\label{full_elements_in_g}
\it Soit
$G=<a,b'><a,b>$. Alors: \newline 
 \hspace*{1cm}Si ($m_{bb'}>2$ ou $m_{ab}=\infty$ ou $m_{ab'}=\infty$), alors
$G$ ne contient pas d'\'el\'ement plein.\newline
 \hspace*{1cm}Si ($m_{bb'}=2$, $m_{ab}<\infty$, $m_{ab'}<\infty$), alors
$G$ a exactement deux \'el\'ements pleins, \`a savoir
\begin{displaymath}
\begin{array}{l}
{\Gamma}_{b',a,b}=\langle m_{ab'}-1,a,b'][b,a,m_{ab}-1 \rangle\\
{\rm et} \\
{\Gamma '}_{b',a,b}=\langle m_{ab'}-1,a,b']a[b,a,m_{ab}-1 \rangle.\\
\end{array}
\end{displaymath}
\end{kinnig_lechel}

 {\bf Preuve : } Rappelons que par d\'efinition de la pl\'enitude, 
l'existence d'un \'el\'ement plein implique que tous les coefficients
de la matrice de Coxeter sont finis. Tout \'el\'ement $g$ de $G$ s'\'ecrit
$xy$ avec $x \in <a,b'>, \ y \in <a,b>$. En posant
$\xi=min(x,xa)$ et $\eta=min(y,ay)$ on voit que
$g$ s'\'ecrit $\xi \varepsilon \eta$ avec
$\varepsilon \in \lbrace e;a \rbrace$. Comme
$\xi \lhd \xi a$, on a $j \leq m_{ab'}-1$ tel
que $\xi=\langle j,a,b']$. De m\^eme,
 on a $k \leq m_{ab}-1$ tel que $\eta=[a,b,k \rangle$.
Si $g$ est plein, $g\geq M_{ab'}$ donc 
$j=m_{ab'}-1$ et de m\^eme $k=m_{ab}-1$. D'o\`u
la proposition. \hfill {Q. E. D.}\newline

\begin{kinnig_lechel}\label{obs_when_a_is_not_reg}
\it Supposons que $m_{ab'} \geq 3$
et que la
restriction de $\phi$ \`a $[e,ab'a]$ est confondue
avec celle de $x \mapsto xa$ 
et que $\beta$ n'est pas a-r\'egulier \`a gauche.
Par la remarque \ref{regularity_criterion} il existe t minimal
avec $\phi([b,a,t\rangle )=[b,a,t+1 \rangle,\ 
\phi([a,b,t+1 \rangle)=[b,a,t+2 \rangle, \
t \leq m_{ab}-3$. Alors 
$ab'[b,a,t \rangle$ est un \'el\'ement
minimal de $W\setminus Q$.
\end{kinnig_lechel}

 {\bf Preuve : } Posons $w=ab'[b,a,t \rangle$. La proposition 
\ref{commutativity_on_ps_suffices} donne:
\begin{displaymath}
\forall x<[b,a,t \rangle, \ 
\phi(b'x)=b'\phi(x), \phi(ab'x)=ab'\phi(x).
\end{displaymath}

 Notons $C=coat([b,a,t \rangle)$
Remarquons que si $t\geq 2$, $|C|=2$ et $\phi$ abaisse
un des \'el\'ements de $C$ tout en  rehaussant l'autre \`a 
$[a,b,t \rangle$ ; alors que si $t=1$ on a $C=\lbrace e\rbrace$,
$\phi(e)=a$. Dans tous les cas, il y a un unique \'el\'ement
rehauss\'e par $\phi$ dans $C$ et cet \'el\'ement est envoy\'e sur
$[a,b,t \rangle$.  

On a alors les calculs suivants:\\
\begin{displaymath}
\begin{array}{c}
coat(w)=\lbrace b'[b,a,t \rangle;[a,b,t+1 \rangle \rbrace \cup ab'C\\
\phi(coat(w))=\lbrace b'[b,a,t+1 \rangle;
[b,a,t+2 \rangle \rbrace \cup ab'\phi(C)\\
Z(\phi,w)=\lbrace w;b'[b,a,t+1 \rangle;[b,a,t+2 \rangle;
ab'[a,b,t \rangle \rbrace\\
\end{array}
\end{displaymath}
 
  Supposons par l'absurde que $w\in Q$. Alors comme
$w\not\in coat(\phi(w))$ (cf. deux\-i\`eme ligne ci-dessus) on a
n\'ecessairement $w\lhd \phi(w)$ donc 
$coat(\phi(w))=Z(\phi,w)$. Soit $m$ un mot r\'eduit repr\'esentant
$\phi(w)$. Comme $[b,a,t+2 \rangle \in {\cal I}$ 
et $b' \leq \phi(w)$, on obtient
$m$ en ins\'erant le caract\`ere $b'$ quelque part dans le mot
$[b,a,t+2 \rangle$.
Maintenant, le fait 
$b'[b,a,t+1 \rangle \lhd \phi(w)$ impose
au $b'$ d'\^etre avant
le premier $a$ apparaissant dans $[b,a,t+1 \rangle$, 
tandis que $w \lhd \phi(w)$ emp\^eche
ce m\^eme $b'$ d'\^etre de ce cot\'e, ce qui est une contradiction
manifeste. \hfill {Q. E. D.}\\

\begin{kinnig_lechel}\label{obs_when_b_is_not_reg}
\it Supposons que $m_{bb'}\geq 3$ et
que la restriction de $\phi$ \`a $[e,b'a]$ est confondue
avec celle de $x \mapsto xa$ et
que $\beta$ n'est pas b-r\'egulier \`a gauche.
Par la remarque \ref{regularity_criterion} il existe t minimal
tel que $\phi([a,b,t \rangle)=[a,b,t+1 \rangle,\ 
\phi([b,a,t+1 \rangle)=[a,b,t+2 \rangle, \
t \leq m_{ab}-3$. Alors $bb'[a,b,t \rangle$ est un \'el\'ement
minimal de $W\setminus Q$.
\end{kinnig_lechel}

 {\bf Preuve : } Posons $w=bb'[a,b,t \rangle$. La proposition 
\ref{commutativity_on_ps_suffices}. donne:
\begin{displaymath}
\forall x<[a,b,t \rangle, \ \phi(b'x)=b'\phi(x), 
\phi(bb'x)=bb'\phi(x).
\end{displaymath}

 Notons $C=coat([a,b,t \rangle)=\lbrace [a,b,t-1 \rangle;
[b,a,t-1 \rangle \rbrace$
(car $t\geq 2$). Remarquons que $\phi$ abaisse
un des \'el\'ements de $C$ et rehausse l'autre \`a 
$[b,a,t \rangle$. On a alors les calculs suivants:\\
\begin{displaymath}
\begin{array}{c}
coat(w)=\lbrace b'[a,b,t \rangle;[b,a,t+1 \rangle \rbrace \cup bb'C\\
\phi(coat(w))=\lbrace b'[a,b,t+1 \rangle;
[a,b,t+2 \rangle \rbrace \cup bb'\phi(C)\\
Z(\phi,w)=\lbrace w;b'[a,b,t+1 \rangle;[a,b,t+2 \rangle;
bb'[b,a,t \rangle \rbrace\\
\end{array}
\end{displaymath}
 
  Supposons par l'absurde que $w\in Q$. Alors comme
$w\not\in \phi(coat(w))$ (cf. deux\-i\`eme ligne ci-dessus) on a
n\'ecessairement $w\lhd \phi(w)$ donc 
$coat(\phi(w))=Z(\phi,w)$. Soit $m$ un mot r\'eduit repr\'esentant
$\phi(w)$. Comme $[a,b,t+2 \rangle \in {\cal I}$ et 
$b' \leq \phi(w)$, on obtient
$m$ en ins\'erant le caract\`ere $b'$ quelque part dans le mot
$[a,b,t+2 \rangle$. Maintenant, le fait 
$w \lhd \phi(w)$ 
impose $\phi(w)=abb'[a,b,t \rangle=abab'[b,a,t-1 \rangle$.
C'est incompatible avec $b'[a,b,t+1 \rangle \lhd \phi(w)$. \hfill {Q. E. D.}\\

 \subsection{\'Etude du cas crois\'e en rang 3. }\label{crossed}

 \hed Dans cette section, on prend $S=\lbrace a;b;b' \rbrace$,
$m_{ab}\geq3,\ m_{ab'}\geq3$, 
et
$\phi(e)=a, \ \phi(b)=ab,\ \phi(b')=b'a$ (cas
``crois\'e''). On note $\beta$ ( ${\beta}^{'}$) 
la restriction de $\phi$ \`a $<a,b>$
(respectivement $<a,b'>$). \newline
 
\setlength{\unitlength}{1cm}

\begin{picture}(12,4)

\put(2.5,2){\circle{.2}}
\put(4.5,2){\circle{.2}}
\put(6.5,2){\circle{.2}}
\put(2.6,2){\line(1,0){1.75}}
\put(4.6,2){\line(1,0){1.75}}
\put(2.5,2.2){$b$}
\put(4.5,2.2){$a$}
\put(6.5,2.2){$b'$}
\put(7.5,3){$\phi(e)=a$}
\put(7.5,2){$\phi(b)=ba$}
\put(7.5,1){$\phi(b')=ab'$}
\put(3.5,0){\bf Figure 1 : Cas crois\'e}

\end{picture}

\bigskip

 Par le th\'eor\`eme \ref{not_too_big}. on a $Q \subseteq <a,b'><a,b>$.
Alors par la proposition \ref{full_elements_in_g}, $Q$ ne pourra \^etre plein
que si $m_{bb'}=2,m_{ab}<\infty,m_{ab'}<\infty$ ce que l'on
suppose dans toute la suite de cette section \ref{crossed}.
De plus, chaque fois que l'on trouve
une nouvelle obstruction $x \not\in Q$ avec $x \leq {\Gamma}_{b',a,b}$
alors on peut en conclure que $\phi$ n'est pas plein. C'est ce que l'on fait  
dans les trois premi\`eres propositions qui suivent.\newline

\begin{kinnig_lechel}\label{roses_part_1}
\it Supposons
$\phi(ab')=ab'a$ (c'est le cas par exemple
si  $m_{ab'}=3$), et que $\beta$ n'est pas a-r\'egulier \`a gauche.
Par la remarque \ref{regularity_criterion} il existe t minimal
avec $\phi([b,a,t\rangle )=[b,a,t+1 \rangle,\ 
\phi([a,b,t+1 \rangle)=[b,a,t+2 \rangle, \
t \leq m_{ab}-3$. Alors 
$ab'[b,a,t \rangle$ est un \'el\'ement
minimal de $W\setminus Q$, donc  $\phi$ n'est pas plein.
\end{kinnig_lechel}

 {\bf Preuve : } C'est la proposition 
\ref{obs_when_a_is_not_reg}. \hfill {Q. E. D.}\newline

\begin{kinnig_lechel}\label{roses_part_2}
\it Supposons
$m_{ab'}\geq 4, 
\phi(ab')=b'ab',\phi(ab'a)=b'ab'a$ (c'est le
cas par exemple si $m_{ab'}=4$) et que $\beta$ 
n'est pas confondu avec l'application $x \mapsto ax$.
Alors il existe t minimal
avec $\phi([b,a,t \rangle)=
[b,a,t+1 \rangle, \
t \leq m_{ab}-2$. 
Dans ces conditions, $ab'[b,a,t \rangle$ est un \'el\'ement
minimal de $W\setminus Q$, donc $\phi$ n'est pas plein.
\end{kinnig_lechel}

 {\bf Preuve : } Posons $w=ab'[b,a,t \rangle$. On a, 
par des calculs analogues
\`a ceux de la proposition \ref{obs_when_a_is_not_reg} :
\begin{displaymath}
\begin{array}{c}
coat(w)=\lbrace b'[b,a,t \rangle;
[a,b,t+1 \rangle;ab'[a,b,t-1 \rangle;ab'[b,a,t-1 \rangle \rbrace \\
\phi(coat(w))=
\lbrace b'[b,a,t+1 \rangle;\phi([a,b,t+1 \rangle );
b'ab'[a,b,t-1 \rangle;b'ab'[b,a,t-1 \rangle 
\rbrace  \\
Z(\phi,w)=
\lbrace w;b'[b,a,t+1 \rangle;\phi([a,b,t+1 \rangle );
b'ab'[a,b,a,t-1 \rangle;b'ab'[b,a,t-1 \rangle 
\rbrace\\

\end{array}
\end{displaymath} 

 On conclut ensuite en disant que si $w\in Q$, comme
$\phi([a,b,t+1 \rangle)$ est dans $\cal I$,
que si $m$ est un mot r\'eduit repr\'esentant $\phi(w)$,
alors $m$ se d\'eduit de $\phi([a,b,t+1 \rangle)$ en ins\'erant un $b'$
quelque part. Mais alors on contredit
$b'ab'[b,a,t-1 \rangle \leq \phi(w)$ 
car $m_{a,b'}\geq 4$.  \hfill {Q. E. D.}\\

\begin{kinnig_lechel}\label{roses_part_3}
\it Supposons
$m_{a,b'}\geq 5, \phi(ab')=b'ab',\phi(ab'a)=ab'ab'$.  
Alors $ab'ba$ est un \'el\'ement
minimal de $W\setminus Q$, donc $\phi$ n'est pas plein.
\end{kinnig_lechel}

 {\bf Preuve : } Posons $w=ab'ba$. On a :
\begin{displaymath}
\begin{array}{l}
coat(w)=\lbrace b'ba;aba;ab'a;ab'b \rbrace \\
\phi(b'ba)=b'\phi(ba)\\
\phi(ab'a)=ab'ab' \ ( {\rm hypoth\grave{e}se } ) \\
\phi(ab'b)=\phi(ab')b=b'ab'b\\
\phi(coat(w))=
\lbrace b'\phi(ba);\phi(aba);ab'ab';b'ab'b \rbrace 
\end{array}
\end{displaymath} 

 Si $m$ est un mot r\'eduit repr\'esentant $\phi(w)$, 
comme $ab'ab' \lhd \phi(w)$ on voit que $m$ s'obtient en
ins\'erant le caract\`ere $b$ dans un mot
$\mu \in \lbrace b'ab'a;ab'ab'\rbrace$. Comme
$b'ab'b \in {\cal I}$ on a 
$\phi(w) \in \lbrace b'ab'ba;b'ab'ab;ab'ab'b\rbrace$. Comme
$w=ab'ba \in {\cal I}$ on a m\^eme 
$\phi(w)=b'ab'ba$.\\ 

 Mais alors on contredit $ab'ab' \leq \phi(w)$ 
car $m_{a,b'}\geq 5$.  \hfill {Q. E. D.}\\

\begin{kinnig_lechel}\label{roses_beta_and_beta_prime_are_reg}
\it Si $\phi$ est plein,
alors $\beta$ est a-r\'egulier \`a gauche
et ${\beta}^{'}$ est a-r\'egulier \`a droite.
\end{kinnig_lechel}

 {\bf Preuve : } Lorsque l'on met bout \`a bout les 
propositions 6.II.1, 6.II.2, et 6.II.3, on s'aper\c{c}oit qu'on a
montr\'e en particulier que si $\phi$ est plein,
alors $\beta$ est $a$-r\'egulier \`a gauche. Par
sym\'etrie, ${\beta}^{'}$ doit aussi \^etre
$a$-r\'egulier \`a droite. \hfill {Q. E. D.}\\  

\begin{kinnig_lechel}\label{roses_beta_or_beta_prime_is_mult}
\it Supposons
que $\beta$ est a-r\'egulier \`a gauche et
que ${\beta}^{'}$ est a-r\'egulier \`a droite; que
$\beta$ n'est pas confondu avec $x\mapsto ax$ et que
${\beta}^{'}$ n'est pas confondu avec $x\mapsto xa$. Il existe
donc $t$ et $t'$ minimaux tels que 
$\phi(\langle t',a,b'])=\langle t'+1,a,b']$ et 
$\phi([b,a,t \rangle)=[b,a,t+1 \rangle$.
 Dans ces conditions, $\langle t',a,b'][b,a,t \rangle$ est un \'el\'ement
minimal de $W\setminus Q$, donc $\phi$ n'est pas plein.
\end{kinnig_lechel}

{\bf Preuve :} Posons $w=\langle t',a,b'][b,a,t \rangle$ 
(remarquons que $t,t' \geq 2$).
On a (en it\'erant la proposition \ref{commutativity_on_ps_suffices} 
pour les quatre derni\`eres
lignes) \\
\begin{displaymath}
\begin{array}{l}
\begin{array}{lcc}
coat(w)=&\lbrace \langle t'-1,a,b'][b,a,t \rangle;&
\langle t'-1,b',a][b,a,t \rangle; \\
 &\langle t',a,b'][b,a,t-1 \rangle; &
\langle t',a,b'][a,b,t-1 \rangle\rbrace 
\end{array}\\
\begin{array}{lllll}
\phi(\langle t'-1,a,b'][b,a,t \rangle)&=&
\langle t'-1,a,b']\phi([b,a,t \rangle)&=&
\langle t'-1,a,b'][b,a,t+1 \rangle, \\ 
\phi(\langle t{'}-1,b',a][b,a,t \rangle)&=&
\langle t{'}-1,b',a]\phi([b,a,t \rangle)&=&
\langle t'-1,b',a][b,a,t+1 \rangle, \\
\phi(\langle t,a,b'][b,a,t-1 \rangle)&=&
\phi(\langle t,a,b'])[b,a,t-1 \rangle&=&
\langle t'+1,a,b'][b,a,t-1 \rangle, \\
\phi(\langle t,a,b'][a,b,t-1 \rangle)&=&
\phi(\langle t,a,b'])[a,b,t-1 \rangle&=&
\langle t'+1,a,b'][a,b,t-1 \rangle.
\end{array}
\end{array}
\end{displaymath}

 Supposons par l'absurde que $w\in Q$. Alors, vu l'\'enumeration
qu'on vient de faire $w$ n'est pas dans $\phi(w)$ donc
$w \lhd \phi(w)$. Comme $w \in {\cal I}$, si $m$ est un
mot r\'eduit repr\'esentant $\phi(w)$, alors $m$ s'obtient
en ins\'erant un certain caract\`ere $c$ dans 
$\langle t',a,b'][b,a,t \rangle$.
Remarquons que $t' \leq m_{ab'}-2$ car 
$\phi(\langle t,a,b']) \neq \langle t,a,b']a$ 
et de m\^eme $t \leq m_{ab}-2$. Alors,
comme $\langle t'+1,a,b']\in {\cal I}$ et 
$[b,a,t+1 \rangle \in {\cal I}$, on
voit que le caract\`ere $c$ devrait \^etre 
\`a la fois au d\'ebut et
\`a la fin de $\langle t',a,b'][b,a,t \rangle$, 
ce qui est absurde. \hfill {Q. E. D.}\\

\begin{kinnig_lechel}\label{roses_equ_for_fullness}
\it Les deux propositions suivantes sont
\'equivalentes :
\begin{displaymath}
\begin{array}{ll}
(1) & \phi {\rm \ est \ plein.}\\
(2) & \left\lbrace
{\begin{array}{l}
\beta  \ {\rm est} \ {\rm r\acute{e}gulier \ \grave{a} \ gauche, \ et \ }
\forall x \in <a,b'> \ \beta^{'}(x)=xa\\
OU \\
{\beta}^{'}  \ {\rm est} \ {\rm r\acute{e}gulier \ \grave{a} \ droite, \ et \ }
\forall y \in <a,b> \ \beta(y)=ay\\
\end{array}}
\right.
\end{array}
\end{displaymath}
\end{kinnig_lechel}

 {\bf Preuve : } En combinant les propositions 
\ref{roses_beta_and_beta_prime_are_reg} et 
\ref{roses_beta_or_beta_prime_is_mult} on voit que
(1) implique (2). R\'eciproquement, par exemple dans la deuxi\`eme
alternative de (2)  on a par le th\'eor\`eme \ref{explicit_formula_for_phi}
\begin{displaymath}
\forall x \in <a,b'>, \ \forall y \in <a,b>, \ 
\phi(xy)=\phi(x)y
\end{displaymath}

 En particulier, on voit que $Q$ contient l'\'el\'ement
${\Gamma}_{b',a,b}$ (cf. proposition \ref{full_elements_in_g})
qui est plein. \hfill {Q. E. D.} \newline

\subsection{\'Etude du cas non d\'egener\'e 
en rang 3. }\label{non_degenerate}

 \hed Dans cette section, on suppose $S=\lbrace a;b;b'\rbrace,
m_{ab}\geq3,\ m_{ab'}\geq3$
(cas
``non d\'egener\'e''). Le cas crois\'e venant
d'\^{e}tre trait\'e \`a la section pr\'ec\'edente, on
suppose ici $\phi(b)=ba, \phi(b')=b'a$. Comme 
pr\'ecedemment, le cas $m_{bb'}>2$ est plus simple.

\begin{picture}(12,4)

\put(2.5,2){\circle{.2}}
\put(4.5,2){\circle{.2}}
\put(6.5,2){\circle{.2}}
\put(2.6,2){\line(1,0){1.75}}
\put(4.6,2){\line(1,0){1.75}}
\put(2.5,2.2){$b$}
\put(4.5,2.2){$a$}
\put(6.5,2.2){$b'$}
\put(7.5,3){$\phi(e)=a$}
\put(7.5,2){$\phi(b)=ba$}
\put(7.5,1){$\phi(b')=b'a$}
\put(2.2,0){\bf Figure 2 : Cas non d\'egener\'e non crois\'e}

\end{picture}

\begin{kinnig_lechel}\label{core_of_nd}
\it Supposons que $\beta$
n'est pas confondu avec $x \mapsto xa$. Alors il
existe $t$ minimal tel que 
$ \phi(\langle t,a,b])=\langle t+1,a,b]$
(avec $m_{ab}\geq t+2 $). 
Soit $H=\lbrace w\in W \ ; \ l(w)=t+1, \ 
b' \leq w, \ \langle t,a,b] \leq w, \
w \neq b'\langle t,a,b] \rbrace$, et $\Gamma$, $\Gamma'$
les \'el\'ements d\'efinis \`a la proposition \ref{full_elements_in_g}. 
Alors \newline
\hspace*{1cm}(1) Si $\phi(ab')=ab'a$, alors 
$\forall w \in H, \ w\not\in Q$.\newline
\hspace*{1cm}(2) Si $\phi(ab')\neq ab'a$, alors 
$abb' \not\in Q, \ ab'b \not\in Q$. \newline
\hspace*{1cm}(3) L'ensemble $Q$ ne contient pas d'\'el\'ement plein, 
sauf dans le cas \nopagebreak $\phi(ab')=ab'a$, $m_{bb'}=2$, $t$ impair. 
Dans ce cas, les seuls (\'eventuels) \'el\'ements 
pleins  \nopagebreak dans $Q$ sont $\Gamma$ et $\Gamma'$.
\end{kinnig_lechel}
  
 {\bf Preuve : } Montrons d'abord (1).

 D\'efinissons un ensemble d'entiers $T$ de la fa\c{c}on
suivante : 
\begin{displaymath}
T=\left\lbrace
\begin{array}{lll}
\hspace*{0cm} [0,t-1] & {\rm si} & m_{bb'}>2 \\
\hspace*{0cm} \lbrace j \in [0,t-2] \ ; \ j \ {\rm pair} \rbrace & 
{\rm si} & m_{bb'}=2  \\
\end{array}
\right.
\end{displaymath} 
 On a alors la description suivante de $H$ : tout
$w\in H$ s'\'ecrit
$w={C}b'{\langle j,a,b]}$, (o\`u $C$ est l'unique
mot tel que l'on ait l'\'egalit\'e de mots
$\langle t,a,b]=C\langle j,a,b]$) 
avec $j \in T$. Soit $v \in coat(w)$.
Si $v=\langle t,a,b]$, on a 
$\phi(v)=\langle t+1,a,b] \neq w$ et
sinon $\phi(v)=va$ (par 
la proposition \ref{equality_on_pp_suffices}, en utilisant le
couplage $x \mapsto xa$ ).
Pour ces derniers $v$ , on a $va \neq w$
(sinon $wa$ serait un coatome de $w$, ce qui contredit
le fait que $w \in \cal I$ et que son dernier caract\`ere
n'est pas un $a$) donc $\phi(v)\neq w $. Ainsi
$w\not\in \phi(coat(w))$, donc
si $w \in Q$ on doit avoir
$w \lhd \phi(w)$. Soit $D$
le mot obtenu en effa\c{c}ant de $C$
son dernier caract\`ere. 
Les \'el\'ements
$x_1=Db'\langle j,a,b]$ et $x_2=\langle t,a,b]$ sont
des coatomes de $w$, et
comme $x_1 \lhd \phi(x_1) (=x_1a)$,
$x_2 \lhd \phi (x_2)$, si $w \in Q$, $y=\phi(w)$
v\'erifie :
\begin{displaymath}
Db'\langle j+1,b,a] \lhd y \ ; \ 
\langle t+1,a,b] \lhd y. 
\end{displaymath}
Remarquons que l'\'el\'ement represent\'e par
$m_2=\langle t+1,a,b]$ est
dans $\cal I$, donc
si $m$ est un mot r\'eduit  
repr\'esentant $y$, alors on obtient $m$
en rajoutant le caract\`ere $b'$ quelque part
dans le mot $m_2$ : 
$m=Ub'V$ et $m_2=UV$. L'\'el\'ement $\phi(x_1)$ 
peut avoir plusieurs \'ecritures r\'eduites
(quand $m_{ab'}=3$) mais parmi celles-ci 
$m_1=Db'\langle j+1,b,a]$ est la seule contenant au plus
un caract\`ere $b'$. Dans tous les cas, $m_1$ doit donc
\^etre un sous-expression de $m$.
On doit donc avoir
$U \geq u=D$ et $V \geq v=\langle j+1,b,a]$.
De plus $V\neq v$ (les derniers caract\`eres diff\`erent)
et du fait des contraintes de longueur 
$U=u$ et $V=vb$. Alors le dernier caract\`ere
de $U$ coinc\i{i}de avec le premier caract\`ere
de $V$, donc $m_2=UV$ n'est pas r\'eduit,
ce qui est absurde. Ceci ach\`eve de montrer $(1)$.\newline 

 Montrons $(2)$. Remarquons que les hypoth\`eses
impliquent que $\phi(ab')=b'ab', \ m_{ab'} \geq 4$.
Soit $w_1=abb'$ ; on a (en utilisant la proposition
\ref{equality_on_pp_suffices}. avec le couplage 
$x \mapsto xa$ pour la derni\`ere
\'egalit\'e ) 
\begin{displaymath}
\begin{array}{l}
coat(abb')=\lbrace ab;ab';bb' \rbrace \\
\phi(ab)\in \lbrace aba;bab \rbrace, \ \phi(ab')=b'ab', \ \phi(bb')=bb'a.
\end{array}
\end{displaymath} 
 Donc si $w_1 \in Q$, on doit avoir $w_1 \lhd \phi(w_1)$ et
pour $m_1$ mot r\'eduit repr\'esentant $\phi(w_1)$ en ins\'erant
un caract\`ere $b$ quelque part dans le mot $b'ab'( \in {\cal I})$.
Alors $m_1$ a au exactement deux 
caract\`eres dans $\lbrace a;b\rbrace$.
C'est incompatible avec $\phi(ab) \lhd \phi(w_1)$. 
Donc $abb' \not\in Q$. Le raisonnement montrant $ab'b \not\in Q$ 
est analogue.\newline

 Montrons maintenant $(3)$. Supposons 
que l'on ait $\zeta \in Q$ plein. \newline

 Cas $\phi(ab')=ab'a$ :\newline\nopagebreak
  
 Soit $p$ le premier caract\`ere de $\langle t,a,b]$.
(ainsi $p=a$ si $t$ est pair et $p=b$ sinon) et
$\bar{p}$ l'autre caract\`ere. Pour $u$ un
pr\'efixe de $\langle t,a,b]$, on note $u^\star$ 
l'unique suffixe de $\langle t,a,b]$ tel que
$\langle t,a,b]=uu^\star$. 
 Commen\c{c}ons par faire la remarque suivante : pour
tout pr\'efixe $u$ de $\langle t,a,b]$, 
\begin{displaymath}
(R)\left\lbrace
\begin{array}{ll}
{\rm Si \ } (m_{bb'}=2,\ t {\rm \ est \ impair),\ alors \ }  &
(ub'u^{\star}\not\in H) \Leftrightarrow (u\in \lbrace e;b \rbrace).\\
{\rm Sinon, \ } &
(ub'u^{\star}\not\in H) \Leftrightarrow (u=e).
\end{array}
\right.
\end{displaymath}

  Consid\'erons une 
occurrence quelconque
de $b'$ dans $\zeta$, i.e. une \'ecriture 
$\zeta=xb'y$, avec $l( \zeta)=l(x)+1+l(y)$.
On sait (cf. la proposition 2.5 de [\ref{du_Cloux_2003}]) que
$[e,x] \cap <a,b>$ 
 a un plus grand
\'el\'ement $\xi$.
 Soit $\xi^{\sharp}$ l'unique \'el\'ement de $<a,b>$
tel que l'on ait $M_{ab}=\xi{\xi^{\sharp}},\ 
m_{ab}=l( \xi)+l( \xi^{\sharp})$. 
Comme $\zeta$ est plein on a $\zeta \geq M_{ab}$.
Par la remarque suivant la proposition 2.5 de [\ref{du_Cloux_2003}],
on en d\'eduit $y \geq {\xi}^{\sharp}$. 
D\'efinissons  un \'el\'ement $\xi'$ de la mani\`ere
suivante : si $p\in D_g(\xi)$, on pose $\xi'=\xi$
et sinon $\xi<M_{ab}$ donc $\xi$ a une
unique \'ecriture r\'eduite. Alors $D_g(\xi)$
est r\'eduit \`a un singleton $\xi_1$ (ou est vide
si $\xi=e$) et on pose alors $\xi'=\xi_1 \xi$ 
(et $\xi'=e$ si $\xi=e$). Posons
aussi $u=min(\langle t,a,b], \xi')$.\newline
 \hspace*{.5cm}Remarquons que $u$ est toujours
un pr\'efixe de $\langle t,a,b]$.
 Montrons maintenant par disjonction ce cas que
$u^{\star} \leq {\xi}^{\sharp}$ :
si $\langle t,a,b] \leq \xi'$, on a 
$u=\langle t,a,b]$ donc $u^{\star}=e$. Si
$\xi' < \langle t,a,b]$ on a
$u=\xi'$ donc 
$l(u^{\star})=t-l(u)=t-l(\xi') \leq t-l(\xi)+1
\leq m_{ab}-1-l(\xi)=l(\xi^{\sharp})-1$
donc $u^{\star}<\xi^{\sharp}$.\newline
 \hspace*{.5cm} Finalement $u \leq \xi \leq x$, 
$u^{\star} \leq {\xi}^{\sharp} \leq y$ donc
$ub'u^{\star} \leq \zeta$ et en particulier
$ub'u^{\star}\not\in H$.\newline
 
\hspace*{.5cm}Vu la remarque plus haut, ceci donne
$u\in \lbrace e;b \rbrace$. Notons que $u=e$ correspond
\`a $\xi \in \lbrace e;\bar{p} \rbrace$ et que $u=b$
correspond \`a $\xi \in \lbrace b;ab \rbrace$.
Distinguons deux cas suivant que $u=e$ pour toute occurrence
de $b'$ dans $\zeta$ ou que l'on puisse avoir $u=b$ parfois.
 Dans le premier cas on a constamment $u=e$, i.e.
$min(\xi',\langle t,a,b])=e$ donc $\xi'=e$ donc
$\xi \in \lbrace e; \bar{p} \rbrace$. En particulier
$\xi\not\geq p$, donc $x \not\geq p$ : derri\`ere
une occurrence de $b'$ dans $\zeta$ on
ne peut avoir de $p$.  Comme $\zeta$ est plein
par rapport \`a $\lbrace p;b' \rbrace$, ceci 
implique $m_{p,b'}=2$. Or, comme $p\in \lbrace a;b \rbrace$
et $m_{a,b'} \geq 3$, ceci implique $p=b$, c'est-\`a-dire
$t$ impair. On a donc : $m_{bb'}=2$,  $t$ impair.
Cette derni\`ere assertion est \'egalement vraie dans le deuxi\`eme
cas o\`u $u$ est parfois \'egal \`a $b$, 
par $(R)$ directement. Donc, dans
tous les cas :\newline
\hspace*{4cm}On a $m_{bb'}=2$, $t$ est impair. \newline

 On peut
\'ecrire $\zeta=xb'y$ avec $y \in <a,b>, \ l(\zeta)=l(x)+1+l(y)$.
On vient de voir que le ``$u$'' de cette d\'ecomposition
est dans $\lbrace e;b \rbrace$, donc le ``$\xi$'' est dans
$\lbrace e;a;b;ab \rbrace$, donc $\xi\not\geq ba$, donc
$x\not\geq ba$.
Comme $\zeta$ est plein on a $x\neq e$ ; on peut alors r\'e\'ecrire $\zeta$
comme $\zeta=x'qb'y$ avec $q \in \lbrace a;b \rbrace ;
l(\zeta)=l(x')+2+l(y)$. On a $q \in \lbrace a;b \rbrace$.
Quitte \`a permuter $b$ et $b'$, on peut supposer $q \neq b$,
donc $q=a$. Alors, comme $x \not\geq ba$ on a $x' \not\geq b$
donc $x' \in <a,b'>$ donc $\zeta \in <a,b'><a,b>$ puis, par
la proposition \ref{full_elements_in_g}, 
$\zeta \in \lbrace \Gamma ; \Gamma' \rbrace$
comme cherch\'e. \newline

 Cas $\phi(ab')\neq ab'a$ :\newline 

Comme $\zeta$ est plein on a $\zeta \geq a$. On a une d\'ecomposition du
type $\zeta=uav$, avec $u \in <b,b'>, \ l(\zeta)=l(u)+1+l(v)$.
Comme $\zeta$ est plein on a $v\neq e$. Alors le premier caract\`ere
$q$ de $v$ est dans $\lbrace b;b' \rbrace$ ; soit $\bar{q}$ 
l'\'el\'ement d\'efini par $\lbrace b;b' \rbrace=\lbrace q;\bar{q} \rbrace$.
On peut \'ecrire $v=qw$ avec $l(v)=1+l(w)$.
Alors comme $\zeta \not\geq aq\bar{q}$ on a $w \not\geq \bar{q}$
donc $w \in <a,q>$. Donc  $\zeta=u(aqw)\in 
<b,b'><a,q>=<q,\bar{q}><a,q>$; comme
$m_{a,\bar{q}} \geq 3$, $\zeta$ ne peut \^etre plein par rapport
\`a $\lbrace a;\bar{q} \rbrace$. \hfill {Q. E. D.}\newline

 En appliquant cette proposition deux fois (la deuxi\`eme fois
en \'echangeant les r\^{o}les de $b'$ et $b$) on voit que 
quand $m_{bb'}>2$, $\phi$
ne peut \^etre plein que si $\beta$ et ${\beta}^{'}$ sont tous deux
des restrictions de $x \mapsto xa$ ; par le th\'eor\`eme 
\ref{mw_is_a_cartesian_product} on
a alors :

\begin{kinnig_lechel}\label{noncommutative_case_for_nd}
\it Si $m_{bb'}>2$, $\phi$ est plein
si et seulement si $\phi$ coincide avec
$x \mapsto xa$.
\end{kinnig_lechel}

 On supposera $m_{bb'}=2$ dans la suite
de cette section.\newline

\smallskip

\begin{kinnig_lechel}\label{one_of_beta_or_beta_prime_is_mult_nd}
\it Si ni $\beta$ ni $\beta^{'}$
ne sont
confondus avec $x \mapsto xa$, alors $\phi$ n'est pas plein.
\end{kinnig_lechel}

 {\bf Preuve : } Supposons par l'absurde que $\phi$ soit plein.
Par la proposition \ref{core_of_nd},
$Q$ contient un unique \'el\'ement plein
de longueur $m_{ab}+m_{ab'}-2$ (que l'on appellera $N$) \`a savoir
$\Gamma_{b',a,b}$. En permutant $b$ et $b'$,
$Q$ contient un unique \'el\'ement plein
de longueur $m_{ab}+m_{ab'}-2$ (que l'on appellera $N'$) \`a savoir
$\Gamma_{b,a,b'}$.
Comme on a $N \neq N^{'}$, (remarquer par exemple
que $bN^{'} \lhd N^{'}$ mais $N \lhd bN$) 
c'est absurde. \hfill {Q. E. D.} \newline

\begin{kinnig_lechel}\label{auxiliary_nd}
\it Supposons que 
$\forall x, \ \beta^{'}(x)=xa$ et que $\beta$ n'est
pas $a$-r\'egulier \`a gauche. Alors $\phi$
n'est pas plein.
\end{kinnig_lechel}
 {\bf Preuve : } La condition ``$\beta$ non $a$-r\'egulier
\`a gauche'' se traduit par  
\begin{displaymath}
\exists t \leq m_{ab}-3, \ \beta([b,a,t \rangle)=[b,a,t+1 \rangle, 
\beta([a,b,t+1 \rangle)=[b,a,t+2 \rangle
\end{displaymath}
 On prend $t$ minimal. La proposition \ref{core_of_nd}
nous dit que si $\phi$ est
plein, alors $Q$ contient
$\Gamma_{b',a,b}$. 

 De plus, la proposition \ref{obs_when_a_is_not_reg}.
montre que $ab'[b,a,t \rangle\not\in Q$ ; comme 
$ab'[b,a,t \rangle \leq \Gamma_{b',a,b}$, c'est impossible car $Q$ est
filtrant d\'ecroissant. \hfill {Q. E. D.}\newline

\begin{kinnig_lechel}\label{equ_for_fullness_nd}
\it Supposons $m_{bb'}=2$,
$\phi$ diff\'erent de $x \mapsto xa$. On a \'equivalence entre
\begin{displaymath}
\begin{array}{ll}
(1)& \phi \ est \ plein\\
(2)& \left \lbrace 
\begin{array}{l}
\grave{A} \ \acute{e}change \ de \ b \ et \ b' \ pr\grave{e}s,\ 
on \ a :\\
 \left \lbrace
\begin{array}{l}
{\phi}_{|<a,b'>} \ est \ simplement \ la \ multiplication \ 
\grave{a} \ droite \ par \ a \\
{\phi}_{|<a,b>} \ est \ a-r\acute{e}gulier \ \grave{a} \
gauche
\end{array}
\right.
\end{array}
\right.
\end{array}
\end{displaymath}
 Le cas \'ech\'eant, $Q$ contient exactement deux \'el\'ements pleins,
\`a savoir $\Gamma_{b',a,b}$ et ${\Gamma'}_{b',a,b}$.
\end{kinnig_lechel}

{\bf Preuve : } L'implication $(2) \Rightarrow (1)$ vient
simplement de la proposition \ref{explicit_formula_for_phi}. 
R\'eciproquement, supposons
$(1)$. Par la proposition 
\ref{one_of_beta_or_beta_prime_is_mult_nd}, $\beta$ ou $\beta^{'}$ est
confondu avec $x \mapsto xa$. Supposons par exemple qu'il
s'agit de $\beta^{'}$. La proposition \ref{auxiliary_nd} assure alors que $\beta$
est $a$-r\'egulier \`a gauche, d'o\`u $(2)$.
 Enfin, si l'un des termes de l'\'equivalence est r\'ealis\'e,
comme $\phi$ est diff\'erent de $x \mapsto xa$, on doit
avoir $\beta$ non confondu avec $x \mapsto xa$, et on peut
d\`es lors utiliser la proposition \ref{core_of_nd}. pour voir
que les seuls (\'eventuels) \'el\'ements pleins de $Q$ sont
$\Gamma$ et $\Gamma'$. Pour v\'erifier qu'ils sont effectivement
dans $Q$, on invoque $<a,b'><a,b> \subseteq Q$, qui provient
de la proposition \ref{explicit_formula_for_phi}.
(avec 
$X=\lbrace a;b' \rbrace, \ 
Y=\lbrace a;b\rbrace$) \hfill {Q. E. D.}\newline

 Il nous reste maintenant \`a examiner le cas d\'egener\'e,
i.e. le cas $\exists i \in \lbrace b;b' \rbrace $,
$m_{ia}=2$. Quitte \`a \'echanger $b'$ et $b$, on peut
supposer $m_{ab'}=2$. Dans ce cas on distingue plusiers
sous-cas suivant les valeurs de $m_{bb'}$.\newline 

 \subsection{\'Etude du cas d\'egener\'e en rang 3. }\label{degenerate}

 \hed Dans cette section, on suppose donc $m_{ab'}=2$. R\`eglons
rapidement le cas $m_{bb'}=2$ par la remarque claire 
suivante:

\begin{merkadenn_lechel}\label{obvious_remark}
\it Si 
Supposons  $m_{ab'}=m_{bb'}=2$.
Dans ce cas on a $W=<a,b> \amalg b'<a,b>$, et tout couplage
est $b'$-r\'egulier \`a gauche, donc d\'efini sur tout
$W$ et plein. 
\end{merkadenn_lechel}

 On suppose donc $m_{bb'} \geq 3$ dans la suite de cette section.

\begin{picture}(12,4)

\put(2.5,2){\circle{.2}}
\put(4.5,2){\circle{.2}}
\put(6.5,2){\circle{.2}}
\put(2.6,2){\line(1,0){1.75}}
\put(4.6,2){\line(1,0){1.75}}
\put(2.5,2.2){$b'$}
\put(4.5,2.2){$b$}
\put(6.5,2.2){$a$}
\put(7.5,2.5){$\phi(e)=a$}
\put(7.5,1.5){$\phi(b)=ba$}
\put(3.5,0){\bf Figure 3 : Cas d\'egener\'e}

\end{picture}

\begin{kinnig_lechel}\label{obs1_d}
\it Supposons que $\beta$
n'est pas confondu avec $x \mapsto xa$. Alors il
existe $t \geq 2$ minimal tel que 
$\phi(\langle t,a,b])=\langle t+1,a,b]$
(donc $m_{ab}\geq t+2$). 
Soit $H=\lbrace w\in W \ ; \ l(w)=t+1, \ 
b' \leq w, \ \langle t,a,b] \leq w, \ w \not\in 
\lbrace b'\langle t,a,b];\langle t,a,b]b' \rbrace \rbrace$.
 Alors 
\begin{displaymath}
\forall w \in H, \ w\not\in Q. \\
\end{displaymath}
\end{kinnig_lechel}
  
 {\bf Preuve : }   Soit $w\in H$. 
Alors, gr\^ace \`a $ab'=b'a$, $w$ s'\'ecrit
$w=Cb'\langle 2j,a,b]$, (o\`u $C$ est 
l'unique mot tel que $\langle t,a,b]=C\langle 2j,a,b]$). 
avec $0 \leq 2j \leq t-1$ (on aurait aussi bien
pu prendre $2j+1$ au lieu de $2j$ dans l'\'ecriture ci-dessus
mais le $2j$ est plus pratique pour la suite). 
 Soit $v \in coat(w)$.
Si $v=\langle t,a,b]$, on a $\phi(v)=\langle t+1,a,b] \neq w$ et
sinon $\phi(v)=va$ (par la proposition \ref{equality_on_pp_suffices}.)
donc $\phi(v)\neq w $ car $a\not\in D_d (w)$. Ainsi
$w\not\in \phi(coat(w))$, donc
si $w \in Q$ on doit avoir
$w \lhd \phi(w)$. 

Les \'el\'ements
$x_1=Cb'\langle 2j-1,a,b]$ et $x_2=\langle t,a,b]$ sont
des coatomes de $w$, et
comme $x_1 \lhd \phi(x_1) (=x_1a)$,
$x_2 \lhd \phi (x_2)$ , si $w \in Q$, $y=\phi(w)$
v\'erifie :
\begin{displaymath}
Cb'\langle 2j,b,a] \lhd y \ ; \ \langle t+1,a,b] \lhd y. 
\end{displaymath}
Consid\'erons les mots
$m_1=Cb'\langle 2j,b,a]$ et $m_2=\langle t+1,a,b]$. 
Si $m$ est un mot r\'eduit  
repr\'esentant $y$, alors on obtient $m$
en rajoutant le caract\`ere $b'$ quelque part
dans le mot $m_2$ : 
$m=Ub'V$ et $m_2=UV$. Maintenant, l'\'el\'ement $\phi(x_1)$
peut avoir plusieurs \'ecritures r\'eduites, mais il n'y en a qu'une 
qui comporte au plus un caract\`ere $b'$, \`a savoir $m_1$.
Donc on doit avoir
$U \geq u=C$ et $V \geq v=\langle 2j,b,a]$.
De plus $V\neq v$ (les derniers caract\`eres diff\`erent)
et du fait des contraintes de longueur, 
$U=u$ et $V=vb$. Alors le dernier caract\`ere de
$U$ coincide avec le premier caract\`ere de
$V$, donc le mot $m_2=UV$ n'est pas r\'eduit,
ce qui est absurde. \hfill {Q. E. D.}\newline

\begin{kinnig_lechel}\label{obs2_d}
\it Supposons $m_{bb'} \geq 4$ et
$\phi(ab)=bab, \ m_{ab} \geq 4$. Alors $abb'b$ est un \'el\'ement
minimal de $W \setminus Q$. 
\end{kinnig_lechel}

 {\bf Preuve : } Posons $w=abb'b$. On a :
\begin{displaymath}
\begin{array}{l}
coat(w)=\lbrace bb'b;abb';ab'b \rbrace \\
\phi(bb'b)=bb'b\phi(e)=bb'ba\\
\phi(abb')=\phi(ab)b'=babb'\\
\phi(ab'b)=\phi(b'ab)=b'\phi(ab)=b'bab
\end{array}
\end{displaymath}

 donc si $w\in Q$ on doit avoir 
\begin{displaymath}
coat( \phi(w))=\lbrace b'bab;abb'b;bb'ba;babb' \rbrace
\end{displaymath}
ce qui est impossible (par exemple il
n'existe pas de $y$ tel que l'on ait \`a la fois 
$b'bab \lhd y, \ abb'b \lhd y $). \hfill{Q. E. D.} \newline
 
\begin{kinnig_lechel}\label{obs3_d}
\it Supposons $m_{bb'}=3$ et
$\phi(ab)=bab, m_{ab} \geq 4$. Alors $abb'ab$ est un \'el\'ement
minimal de $W \setminus Q$. 
\end{kinnig_lechel}

 {\bf Preuve : } Posons $w=abb'ab$. On a :
\begin{displaymath}
\begin{array}{l}
coat(w)=\lbrace bb'ab;abab;abb'b;abb'a \rbrace \\
\phi(bb'ab)=bb'\phi(ab)=bb'bab \\
\phi(abb'b)=\phi(b'abb')=b'\phi(ab)b'=b'babb'\\
\phi(abb'a)=\phi(abab')=\phi(aba)b' 
\end{array}
\end{displaymath}

 donc si $w\in Q$ on doit avoir $w \lhd \phi(w)=y$ 
et $b'babb' \lhd y, \ bb'bab \lhd y$. 
Comme $b'babb' \in {\cal I}$, 
et que $bb'bab$ a exactement trois \'ecritures
r\'eduites \`a savoir $bb'bab, \ b'bb'ab$ et $b'bab'b$,
on en d\'eduit 
$y=b'bab'bb'$ qui est incompatible avec
$\phi(aba)b'\lhd y$. \hfill{Q. E. D.} \newline

\section{\bf Cas g\'en\'eral. }\label{gen}

\hed  On consid\`ere maintenant un couplage maximal $(Q,\phi)$ sur
un syst\`eme de Coxeter $(W,S)$ quelconque. On va progressivement
montrer que $\phi$ est r\'eductible dans tous les
cas. Bien entendu, on peut supposer
que $\phi$ est plein. Par le lemme \ref{subgroups_remain_full}, si $<J>$
est un sous-groupe parabolique stable par $\phi$, alors
${\phi}_{|<J>\cap Q}$ reste plein, ce qui va nous
permettre d'utiliser les r\'esultats d\'eja obtenus
en rang 3. Posons $a=\phi(e)$,

\begin{displaymath}
\begin{array}{l}
C=\lbrace s \in S \ ; \ sa=as \rbrace \\
U=\lbrace s \in S \setminus C \ ; \ \phi(s)=sa \rbrace \\
V=\lbrace s \in S \setminus C \ ; \ \phi(s)=as \rbrace \\
\end{array}
\end{displaymath}

 Commen\c{c}ons par traiter le cas dit ``crois\'e'' : \newline

\begin{kinnig}\label{crossed_case}
\it Supposons que $\phi$ est un couplage
``crois\'e'' (i.e. tel que $U\neq \emptyset, V\neq\emptyset$).  
Alors $\phi$ est r\'eductible. Plus pr\'ecisement, quitte
\`a \'echanger $U$ et $V$, quitte \`a \'echanger 
la gauche et la droite,
on a :
\begin{displaymath}
\begin{array}{l}
Q=<U\cup C>(<V \cup C> \cap Q) \\
\forall (x,y) \in <U \cup C> \times (<V \cup C> \cap Q), \
\phi(xy)=x\phi(y).
\end{array}
\end{displaymath}
\end{kinnig}

 {\bf Preuve : } Disons que $u \in U$ est
``inerte'' si ${\phi}_{<a,u>}$ est $a$-r\'egulier
\`a gauche, et que $u$ est
``fortement inerte'' si ${\phi}_{<a,u>}$ est simplement
la restriction de $x \mapsto xa$ \`a $<a,u>$.
De fa\c{c}on sym\'etrique, 
disons que $v \in V$ est
inerte si ${\phi}_{<a,v>}$ est $a$-r\'egulier
\`a droite, et que $v$ est
fortement inerte si ${\phi}_{<a,v>}$ est simplement
la restriction de $x \mapsto ax$ \`a $<a,v>$.

 Par la proposition \ref{roses_beta_or_beta_prime_is_mult}, on voit que

\begin{displaymath}
\forall (u,v) \in U \times V, 
\left\lbrace
\begin{array}{l}
( u \ { \rm est \ inerte}, \ 
v \ { \rm est \ fortement \ inerte})\\ 
{\rm ou} \\
( u \ { \rm est \ fortement \ inerte}, \ 
v \ { \rm est \ inerte})\\ 
\end{array}
\right.
\end{displaymath}

donc tous les \'el\'ements de $U\cup V$ sont inertes, et
\begin{displaymath}
\forall (u,v) \in U \times V, \
u {\rm \ ou \ } v { \rm \ est \ fortement \ inerte}
\end{displaymath}

d'o\`u l'on d\'eduit ais\'ement que l'un des 
deux ensembles $U$, $V$ n'est compos\'e
que d'\'el\'ements fortement inertes. Supposons
par exemple que ce soit $U$. 

 Alors les th\'eor\`emes \ref{explicit_formula_for_phi} 
et \ref{not_too_big}. donnent une
\'egalit\'e par double inclusion pour
$Q$ : le th\'eor\`eme \ref{not_too_big}. donne
$Q \subseteq <U \cup C><V \cup C>$, donc
$Q \subseteq <U \cup C>(<V \cup C> \cap Q)$ car $Q$ est
d\'ecroissant, et le th\'eor\`eme \ref{explicit_formula_for_phi} donne
$ <U \cup C>(<V \cup C> \cap Q) \subseteq Q$. 

 Expliquons maintenant pourquoi ceci implique que
$\phi$ est r\'eductible : soit $\omega$ une orbite
pleine, $\omega=\lbrace m;M \rbrace$ avec $M=\phi(m)$
et $M$ plein. Alors il existe 
$(x,y)\in (U \cup C) \times ((V \cup C) \cap Q)$ tels
que $m=xy, \ M=x\phi(y)$. On peut
supposer $l(m)=l(x)+l(y)$ par la
r\`egle de l'effacement. Il est facile de voir
que pour toute partie $J$ de $S$ contenant $a$,
$<J> \cap Q$ est stable par $\phi$. En particulier
$\phi(y) \in <(V cup C) \cap Q>$. Comme $U\neq \emptyset$
et $M$ est plein, on en d\'eduit $x\neq e$. Soit 
$x_1 \in D_g(x)$ ; alors $x_1$ est r\'egulier \`a gauche
(car $x_1 \in U \cup C$) et $x_1$ est dans l'ensemble
de descente \`a gauche de $m$, donc l'orbite
$\omega$ est r\'eductible. \hfill {Q. E. D.}\newline

\begin{heuliadenn}\label{simply_laced_case}
\it Tout couplage d\'efini
sur un groupe de Coxeter simplement enlac\'e
est r\'eductible.
\end{heuliadenn}

 {\bf Preuve : } Soit $\phi$ un tel couplage ; on
peut prendre $\phi$ maximal est plein. Si $\phi$ est
crois\'e, alors on utilise la proposition qui 
pr\'ec\`ede. Sinon, l'hypoth\`ese sur le
syst\`eme de Coxeter implique que 
$\forall s \in S \setminus \lbrace a \rbrace, \ 
\forall x \in <a,s>, \ \phi(x)=xa$,  donc 
par le th\'eor\`eme \ref{mw_is_a_cartesian_product}
$\forall x \in W, \ \phi(x)=xa$
auquel cas le r\'esultat est clair. \hfill {Q. E. D.}\newline

 Revenant au cas g\'en\'eral, on voit que l'on peut
toujours se ramener au cas  ou par exemple $V=\emptyset$,
i.e. 
\begin{displaymath}
\forall s \in S, \phi(s)=sa.  
\end{displaymath}

En utilisant la proposition \ref{one_of_beta_or_beta_prime_is_mult_nd}, 
on peut m\^eme supposer
que pour tout $s\in S$ except\'e au plus un \'el\'ement,
\begin{displaymath}
\forall x \in <s,a>, \ \phi(x)=xa.
\end{displaymath}

 Bien s\^ur, le cas non-trivial est le cas o\`u il
existe effectivement un \'el\'ement (que l'on
notera $b$ ) tel que
$\exists x \in <a,b>, \ \phi(x) \neq xa$. Changeant l\'eg\`erement
de notation afin de travailler avec des sous-ensembles
disjoints de $S$, posons

\begin{displaymath}
\begin{array}{l}
C= \lbrace s \in S \ ; \ sa=as, \ s\neq a \rbrace \\
U= S \setminus (C \cup \lbrace a;b \rbrace)
\end{array}
\end{displaymath}

On est alors dans la situation suivante :

\begin{displaymath}
\begin{array}{l}
S=C \amalg U \amalg \lbrace a;b \rbrace \\
\forall c \in C,\ ca=ac, \ c \ est \ r\acute{e}gulier \ 
\grave{a} \ gauche \ et \ \grave{a} \ droite \\
\forall u \in U,\ m_{au} \geq 3, ub=bu, \ u \ est \ r\acute{e}gulier \  
\grave{a} \ gauche.
\end{array}
\end{displaymath}
situation dont une bonne partie est r\'esum\'ee par le dessin suivant :

\begin{picture}(9,5)

\put(4.5,2){\circle{.2}}
\put(6.5,2){\circle{.2}}
\put(4.6,2){\line(1,0){1.8}}
\put(4.15,2){$a$}
\put(6.7,2){$b$}
\put(4.5,4){\circle{.7}}
\put(6.5,4){\circle{.7}}
\put(4.4,3.9){$U$}
\put(6.4,3.9){$C$}
\put(4.5,2.13){\line(0,1){1.53}}
\put(6.5,2.13){\line(0,1){1.53}}
\multiput(4.85,4)(0.2,0){7}{\line(1,0){0.1}}
\put(5,1){\bf Figure 4}

\end{picture}

 (pour le ``$ub=bu$'' \`a la derni\`ere ligne, 
utiliser la proposition \ref{noncommutative_case_for_nd} ). Notons 
$C^{'}=\lbrace c\in C \ ; \ m_{bc} \geq 3 \rbrace $,
$C^{''}=C \setminus C'=\lbrace c\in S \ ; \ sa=as, \ sb=bs \rbrace$
 et pour
$t$ entier,
\begin{displaymath}
T_{t}=\langle t,a,b]
\end{displaymath} 

 Comme $\phi$ n'est pas confondu avec
$x \mapsto xa$ il existe $t$ minimal
tel que $\phi(T_t) \neq T_{t}a$,
donc $t \leq m_{ab}-2$ et
$\phi(T_t)=T_{t+1}$. Notons, pour $x \in U \cup C'$,
\begin{displaymath}
H_x=
\left\lbrace
\begin{array}{ll}
\lbrace w \ ; \ T_{t} \lhd w, \ u\leq w, \ 
w\neq uT_{t} \rbrace
& si \  x=u \in U \\
\lbrace w \ ; \ T_{t} \lhd w, \ c'\leq w, \ 
w\not\in \lbrace c'T_{t};T_{t}c' \rbrace \rbrace
& si \  x=c' \in C',t>2, \\
\lbrace abcb \rbrace
& si \  x=c \in C',t=2, m_{bc}>3 \\
\lbrace abcab \rbrace
& si \  x=c \in C',t=2, m_{bc}=3. 
\end{array}
\right.
\end{displaymath}

On a alors,
par les propositions  \ref{core_of_nd}, \ref{obs1_d}, \ref{obs2_d} et 
\ref{obs3_d} : \newline\nopagebreak

 \hspace*{1cm} Si $t\geq 2$, on a $\forall x \in U\cup C', \  H_{x} \cap Q= 
\emptyset$ \hfill (1) \newline

 Notons $H$ la r\'eunion des diff\'erents $H_x$.
Supposons par l'absurde que $\phi$ n'est
pas r\'eductible. 
Alors  $|S|>2$ et il existe une orbite pleine non
r\'eductible, i.e. il existe $\zeta \in Q$ avec
$\zeta \lhd \phi( \zeta)$, $\phi( \zeta)$ plein
tel que $\zeta$ soit {\bf irr\'eductible}, i.e. 
tel que 
$D_{g}( \zeta )$ ne contienne pas
d'\'el\'ement r\'egulier \`a gauche et
 $D_{d}( \zeta )$ ne contienne pas
d'\'el\'ement r\'egulier \`a droite. Alors
nous affirmons que $\zeta \geq \langle m_{ab}-2,a,b]$.

 En effet, consid\'erons 
l'ensemble 
$A=\lbrace y \leq \zeta \ ; \ 
y\lhd \phi(y), \ \phi(y) \geq M_{ab}\rbrace$.
L'ensemble $A$ est non vide car il contient $\zeta$.
Soit $y_0$ un \'el\'ement minimal de $A$.
On a $\phi(y_0) \geq M_{ab}$. Si $\phi(y_0)=M_{ab}$,
comme le sous-groupe dih\'edral principal $<a,b>$ 
est stable par $\phi$, $y_0$ est un \'el\'ement de
$<a,b>$ de longueur $m_{ab}-1$, donc
$y_0 \geq \langle m_{ab}-2,a,b]$ puis 
$\zeta  \geq \langle m_{ab}-2,a,b]$.
Sinon on a $\phi(y_0)>M_{ab}$ et il existe
alors $z$ tel que $M_{ab} \leq z \lhd \phi(y_0)$.
On sait que $z \in Z(\phi,y_0)=
\lbrace y_0 \rbrace \cup
\phi(B)$
o\`u l'on a pos\'e 
$B=\lbrace x \lhd y_0 \ ; \ x \lhd \phi(x) \rbrace$.
Or si $z=\phi(y_1)$ avec $y_1 \in B$, on aurait
$y_1 \in A$ ce qui contredirait la minimalit\'e de $y_0$.
Donc $z=y_0$. Alors $M_{ab} \leq y_0 \leq \zeta$  
ce qui donne encore $\zeta \geq \langle m_{ab}-2,a,b]$.\newline

 Remarquons que le sous-groupe $G=<\lbrace a,b \rbrace \cup C''>$ 
de $W$ est isomorphe au produit commutatif de groupes
$<a,b> \times <C''>$, et que de plus on a $\forall(x,y) \in <a,b> \times C'',\
\phi(xy)=\phi(x)y$ donc la restriction de $\phi$ \`a $G$ est r\'eductible.
En particulier $\zeta\not\in G$.
Ainsi $\zeta$ poss\`ede les propri\'et\'es suivantes :\newline

\hspace*{1.5cm}
$( \zeta \geq \langle m_{ab}-2,a,b], \
\zeta \ {\rm irr\acute{e}ductible}, \ 
\zeta \not\in <\lbrace a;b \rbrace \cup C''>)$
\hfill (3) \newline

 Soit maintenant
$k$ un param\`etre entier tel que $t \leq k \leq m_{ab}-2$
et $k \equiv t (mod \ 2)$ (on peut toujours 
trouver des $k$ v\'erifiant l'ensemble de
ces conditions, par exemple $k=t$ ). 
On a une d\'ecomposition 
\begin{displaymath}
\begin{array}{l}
\hspace*{2cm}\left\lbrace
\begin{array}{l}
\zeta=ud_1x_1d_2x_2 \ldots d_{k-1}x_{k-1}d_kv, \\
l(\zeta)=l(u)+\sum_{i}{l(x_i)}+l(v)+k \\
d_1d_2 \ldots d_k=\langle k,a,b]
{\rm  \ ( \acute{e}galit\acute{e} \ 
caract\grave{e}re \ \grave{a} \ caract\grave{e}re).}
\end{array}
\right\rbrace \hspace*{.7cm} (4)
\end{array}
\end{displaymath}

 Nous allons montrer qu'on 
peut en fait supposer $\forall i\in [1,k-1], \ x_i=e$
dans cette d\'ecomposition.
Montrons ceci par r\'ecurrence sur $i$\nolinebreak[2] ; pour
chaque $i$ on consid\`ere une \'ecriture de
la forme (4) dans laquelle on a d\'eja $\forall j<i, \ x_j=e$
(hypoth\`ese de r\'ecurrence) et dans laquelle aussi
$l(x_i)$ est minimal. On suppose par l'absurde
$x_i \neq e$ et on prend $s \in D_g(x_i)$.  \newline

 Remarquons d'abord que $s\notin \lbrace a;b \rbrace$ car
sinon comme l'\'ecriture $d_ix_i$ est r\'eduite on a n\'ec\'essairement
$s=d_{i+1}$, et en rempla\c{c}ant $(d_i,x_i,d_{i+1},x_{i+1})$ par
$(d_i,e,d_{i+1},(d_{i+1}x_i)d_{i+1}x_{i+1})$ (\'etant compris
que $x_{i+1}$ s'appelle $v$ si $i=k-1$) on contredit la minimalit\'e
de $l(x_i)$. \newline

 Supposons $s\in U$. Alors 
$t$ est impair, par
le corollaire 6.III.2. En particulier $d_1=b$
car $k \equiv t (mod \ 2)$. Si 
$i=1$ on contredit la minimalit\'e de
$x_1$ en rempla\c{c}ant $(u,x_1)$ par
$(us,sx_1)$.
Si $i>1$ alors 
$\zeta \geq h=d_1 \ldots d_i s d_{i+1} \ldots d_t \in H$ contredit
(1) (cas extr\^eme : si $i \geq t$ on a 
$h=d_1 \ldots d_ts=\langle t,a,b]s$).\newline

 Supposons $s \in C'$. 
Si $i>1$ (ce cas n'existant que
quand $t\geq 3$) alors 
en posant 
$h=d_1 \ldots d_i s d_{i+1} \ldots d_t$ si $i \leq t-1$
et $h=d_1 \ldots d_{t-1}sd_k$ si $i \geq t$, 
\`a chaque fois $\zeta \geq h \in H$ contredit (1).
Reste le cas $i=1$. Si $d_1=a$, on obtient une
contradiction sur la minimalit\'e de $x_1$ en
rempla\c{c}ant $(u,x_1)$ par $(us,sx_1)$. Sinon
$d_1=b$, donc $t$ et $k$ sont impairs, ( et
en particulier $t \geq 3$ ) alors en posant
$h=d_1sd_2 \ldots d_t$ la relation
$\zeta \geq h\in H$ contredit (1).\newline
 
 Enfin, si $s \in C \setminus C'$, ce qui veut
dire que $s$ commute \`a la fois avec $a$ et $b$,
alors en rempla\c{c}ant $(u,x_i)$ par $(us,sx_i)$
on contredit la minimalit\'e de $x_i$.\newline

 Dans tous les cas on est arriv\'e \`a obtenir une
contradiction ; c'est donc que $\forall i, \ x_i=e$.
On a alors
\begin{displaymath}
\begin{array}{l}
\hspace*{4cm}\left\lbrace
\begin{array}{l}
\zeta=u\langle k,a,b]v, \\
l(\zeta)=l(u)+k+l(v) 
\end{array}
\right\rbrace \hspace*{3.3cm} (5)
\end{array}
\end{displaymath}

 Nous affirmons que $v \in <a,b>$. En effet, si $v=e$
c'est vrai ; sinon prenons $s \in D_d(v)$ ; on
peut \'ecrire $v=v's$ avec $v'<v$. On a $s\not\in U$
car $\forall y \in U, \ \langle t,a,b]y \in H$. De
plus $s\not\in C$ car $D_d(\zeta)$ ne contient pas d'\'el\'ements 
r\'eguliers \`a droite. Donc $s\in \lbrace a;b \rbrace$.
Supposons $s=b$. 
Alors 
$\supp(v') \cap (U \cup C')= \emptyset$ 
( sinon, si $c \in  \supp(v') \cap (U \cup C')$, on a
$\zeta \geq d_1 \ldots d_{t-1}cb  \in H$
qui est exclu par $(1)$). Donc
$\supp(v') \subseteq \lbrace a;b \rbrace \cup C''$ ; donc
$v'$ s'\'ecrit $v'=v''c$ avec $l(v)=l(v'')+l(c), \ 
v'' \in <a,b>, \ c \in <C''>$. On a alors $\zeta=u\langle t,a,b]v''bc$.
Comme $D_d(\zeta)$ ne contient pas d'\'el\'ements 
r\'eguliers \`a droite, cela donne $c=e$ donc $v=v''b \in <a,b>$. 
Reste \`a traiter le cas $s=a$. Si $v=a$ on a encore
$v \in <a,b>$. Sinon on prend $s_2 \in D_d(v')$ et on
note $v''=v's_2$. Alors $s_2 \not\in U$ (car 
$\forall y \in U, \ \langle t,a,b]y \in H$) et
$s_2 \in C$ donnerait $\zeta=u\langle t,a,b]v''as_2$, d'o\`u
un \'el\'ement r\'egulier \`a droite dans $D_d(\zeta)$, ce qui
est exclu. Donc $s_2 \in \lbrace a;b \rbrace$, et comme le mot
$s_2s=s_2a$ est r\'eduit on a $s_2=b$. En refaisant alors le raisonnement
fait pour le cas $s=b$ d\'ecal\'e d'un indice, on finit de
v\'erifier que  $v\in <a,b>$ dans tous les cas.\newline

 Comme $\zeta \not\in <\lbrace a,b \rbrace \cup C''>$
il existe $s\not\in \lbrace a;b \rbrace \cup C''$ dans le
support de $u$ \nolinebreak[1]: on a donc 
$s \in U \cup C'$. Soit $u_1 \in D_g(u)$
et $u'=u_1u$. Comme $u_1$ n'est pas
r\'egulier \`a gauche on a $u_1 \in \lbrace a;b \rbrace$.\newline

 Supposons $t$ impair et $u_1=a$. Alors
$a$ n'est pas r\'egulier \`a gauche, 
ce qui implique d\'eja $U=\emptyset$ par
la proposition \ref{core_of_nd}. Donc $s\in C'$.
De plus, si $u_2 \in C$ alors $D_g(\zeta)$ 
contiendrait un \'el\'ement r\'egulier \`a gauche \`a
savoir $u_2$ : impossible. Donc
$u_2 \in \lbrace a;b \rbrace$ puis
$u_2=b$. Alors 
$\zeta \geq bsd_2d_3 \ldots d_t \in H$ ce qui contredit (1).\newline

  Supposons $t$ impair et $u_1=b$. Comme on sait
que $\forall y \in C', \ byd_2 \ldots d_t \not\in H$,
on voit que $\supp(u) \cap C'=\emptyset$, et en particulier
$s\in U$, $u_2\not\in C'$. Si $u_2 \in U\cup C''$ alors $D_g(\zeta)$ 
contiendrait un \'el\'ement r\'egulier \`a gauche \`a
savoir $u_2$ : impossible. Donc
$u_2 \in \lbrace a;b \rbrace$ puis
$u_2=a$. Alors 
$\zeta \geq basd_2d_3 \ldots d_t \in H$ ce qui contredit (1).\newline

 Supposons $t$ pair  et $u_1=b$. Alors
$b$ n'est pas r\'egulier \`a gauche, donc
on a $r \leq m_{ab}-3$ minimal tel que
$\phi([a,b,r \rangle)=\phi([a,b,r+1 \rangle)$
et $\phi([b,a,r+1 \rangle)=\phi([a,b,r+2 \rangle)$.
Comme $t$ est pair, on a $r\geq t+2$.
Posons $k=r$ si $r$ est pair et 
$k=r+1$ sinon ;
alors $k$ v\'erifie les hypoth\`eses pr\'ecis\'ees au d\'ebut
de cette d\'emonstration, et de plus
$k>r$ donc $\langle k,a,b]>[a,b,r \rangle$
puis $\zeta \geq bu'[a,b,r \rangle$.
Par la proposition \ref{obs_when_b_is_not_reg}, 
$\forall y \in C', \ by[a,b,r \rangle \not\in Q$ ; on
en d\'eduit $\supp(u') \cap C'=\emptyset$, et en particulier
$s\in U, u_2 \not\in C'$. Si $u_2 \in U$ alors $D_g(\zeta)$ 
contiendrait un \'el\'ement r\'egulier \`a gauche \`a
savoir $u_2$ : impossible. Finalement 
$u_2 \in \lbrace a;b \rbrace$, puis  $u_2=a$.
 Alors 
$\zeta \geq asd_2d_3 \ldots d_t \in H$ ce qui contredit (1).\newline

 Supposons $t$ pair  et $u_1=a$. Alors
$a$ n'est pas r\'egulier \`a gauche, ce qui
donne $U=\emptyset$ par la proposition 
\ref{core_of_nd}. Donc $s\in C'$.
Si $u_2 \in C$ alors $D_g(\zeta)$ 
contiendrait un \'el\'ement r\'egulier \`a gauche \`a
savoir $u_2$ : impossible. Donc $u_2 \in \lbrace a;b \rbrace$
puis $u_2=b$.  Alors en posant
$h=absd_2$ si $t=2$ et $m_{bs}>3$,
$h=absd_1d_2$ si $t=2$ et $m_{bs}=3$,
et $h=absd_1d_2 \ldots d_{t-2}$ si $t>2$, on a
dans tous les cas
$\zeta \geq h\in H$ ce qui contredit (1).\newline

 Donc, au bout du compte : \newline

\begin{teorem}\label{the_big_result}
\it Pour tout syst\`eme de
Coxeter (W,S), tout couplage distingu\'e sur
$W$ est r\'eductible.   
\end{teorem}

 En combinant avec la proposition \ref{descent_formulas}, on
obtient imm\'ediatement :

\begin{heuliadenn}\label{final_result}
\it Soit $(W,S)$ un syst\`eme de Coxeter, $\phi$ un couplage
distingu\'e de $W$, $x,y \in W$
tels que $x\lhd \phi(x), \ y \lhd \phi(y)$. Alors \newline
\hspace*{4cm}{$R_{\phi(x),\phi(y)}=R_{x,y}$} \hfill (\ref{final_result}.1)\\
\hspace*{3.5cm}{$R_{x,\phi(y)}=(q-1)R_{x,y}+qR_{\phi(x),y}$}\hfill 
(\ref{final_result}.2)\\ 
\end{heuliadenn}

 Bien que nous n'en ayons pas eu besoin ici, il est int\'eressant
de faire la remarque suivante (on note
${\cal M}(W)$ l'ensemble des couplages maximaux
d'un groupe de Coxeter $W$ et pour $a \in S$, 
${\cal M}_a (W)=\lbrace \phi \in {\cal M}(W) \ 
; \ \phi(e)=a\rbrace$) :\newline

\begin{kinnig}\label{final_remark}
 Soit $(W,S)$ un syst\`eme de Coxeter et $a \in S$.
Alors les seuls \'el\'ements
de ${\cal M}_a(W)$ d\'efinis sur tout $W$ sont les multiplications
\`a gauche et \`a droite par un $a$, sauf dans le cas 
d\'egener\'e 
\begin{displaymath}
S=\lbrace a;b \rbrace \amalg C, \ 
\forall c \in C, m_{ac}=m_{bc}=2.
\end{displaymath}
 Dans ce cas, $W$ est isomorphe au produit commutatif de
groupes de Coxeter 
$<C> \times <a,b>$, tous les \'el\'ements de ${\cal M}_a(W)$
sont $<C>$-r\'eguliers (i.e. v\'erifient
$\forall c \in <C>, \ \forall x \in Q, \ cx \in Q, \ \phi(cx)=c\phi(x)$ )
donc d\'efinis sur tout $W$. De plus, l'op\'eration de restriction
sur $<a,b>$ r\'ealise une bijection de ${\cal M}_a(W)$ sur
${\cal M}_a(<a,b>)$. 
\end{kinnig}

 {\bf Preuve : } Soit $\phi \in {\cal M}_a(W)$.
Posons $C=\lbrace s \in S \ ; \ sa=as \rbrace$, 
$U=\lbrace s\in S\setminus C \ ; \ \phi(s)=sa \rbrace$,
$V=\lbrace s\in S\setminus C \ ; \ \phi(s)=as \rbrace$.
Dans le cas ``crois\'e'', i.e. 
$U \neq \emptyset, \ V \neq \emptyset$, alors 
par la proposition 7.1. $Q \subseteq <U \cup C><V \cup C>$,
donc pour $u \in U, \ v \in V$ on a $vau \not\in Q$ donc
$Q \neq W$. On peut donc supposer $V=\emptyset$,
i.e. $\forall s \in S, \phi(s)=sa$. Si $\phi$ n'est
pas confondu
avec le couplage $c_a$ de
 multiplication \`a droite par $a$, par le
th\'eor\`eme
\ref{mw_is_a_cartesian_product} on a 
$b\in S \setminus \lbrace a \rbrace$ tel que
$\phi_{<a,b>}$ ne soit pas une restriction de
$c_a$, donc tel que $\exists t \leq m_{ab}-2$,
$ \phi( \langle t,a,b])=\langle t+1,a,b] )$. Alors
si $U \neq \emptyset$, on a pour $u\in U$,
$\langle t,a,b]u \not\in Q$ par la proposition \ref{core_of_nd}.
On peut donc supposer $U=\emptyset$. Ainsi
$S \subseteq \lbrace a;b \rbrace \cup C''$ avec
$C''=C \setminus \lbrace a \rbrace$. Notons
$C'=\lbrace c\in C'' \ ; \ c \neq b, \ m_{bc}>2 \rbrace$. Supposons
$C' \neq \emptyset$, et prenons $c \in C'$.\newline

Si $t>2$, on a $\langle t-1,b,a]cb \not\in Q$. 
Si $t=2$ et $m_{bc}>3$ on a $abcb \not\in Q$ par
la proposition \ref{obs2_d}. 
Si $t=2$ et $m_{bc}=3$ on a $abcab \not\in Q$ par
la proposition \ref{obs3_d}. Cette disjonction de cas
montre que $Q \neq W$ d\`es que $C' \neq \emptyset$.\newline

 On peut donc supposer $C'=\emptyset$, c'est-\`a-dire
qu'on est dans le cas d\'egener\'e mentionn\'e par
l'\'enonc\'e. Le reste de la proposition est
clair. \hfill {Q. E. D.}\newline

\bigskip
\bigskip
\bigskip
\bigskip

\centerline{\large Remerciements}
\smallskip
 
Je remercie chaleureusement mon directeur de th\`ese Fokko
du Cloux qui m'a aid\'e tout au long de la pr\'eparation de
cet article et qui s'est \'egalement charg\'e du travail
de relecture et de correction des versions pr\'eliminaires
de ce texte.  
 
\bigskip

\hspace*{4cm}{\large Bibliographie}\\*

\refstepcounter{biblio}

\begin{trivlist}

\item[] [\thebiblio]\label{Brenti_1994} 
 F. Brenti. A combinatorial formula for
Kazhdan-Lusztig polynomials, 
{\it Invent. Math.}, {\bf 118} : 371-394, 1994.
\refstepcounter{biblio}

\item[] [\thebiblio]\label{Brenti_2004} 
 F. Brenti. The intersection cohomology of Schubert
  varieties is a combinatorial invariant, 
{\it Europ. J. Combin.}, in press.
\refstepcounter{biblio}

\item[] [\thebiblio]\label{du_Cloux_2000} 
F. du Cloux. An abstract model for Bruhat intervals.  {\it
    Europ. J. Combin.}, {\bf 21} : 197-222, 2000.
\refstepcounter{biblio}

\item[] [\thebiblio]\label{du_Cloux_Coxeter} 
F. du Cloux. {\tt Coxeter}, version beta. Disponible sur
\url{http://www.desargues.univ-lyon1.fr/home/ducloux/coxeter.html}
\refstepcounter{biblio}

\item[] [\thebiblio]\label{du_Cloux_2003} 
F. du Cloux. Rigidity of Schubert closures and invariance
of Kazhdan-Lusztig polynomials.  {\it
    Adv. in Math. }, {\bf 180} : 197-222, 2003.
\refstepcounter{biblio}

\item[] [\thebiblio]\label{du_Cloux_1999} 
F. du Cloux. A transducer approach
to Coxeter groups, {\it J. of Symb. Comp}, {\bf 27} : 1-14, 1999.
\refstepcounter{biblio}

\item[] [\thebiblio]\label{Dyer_1987} 
M. Dyer. Hecke Algebras and reflections in Coxeter groups.
  PhD thesis, University of Sydney, 1987.
\refstepcounter{biblio}

\item[] [\thebiblio]\label{Humphreys} 
J. E. Humphreys. Reflection Groups and Coxeter Groups,
Cambridge Studies in Advanced Mathematics, vol. 29,
Cambridge University Press, Cambridge, 1990.
\refstepcounter{biblio}

\item[] [\thebiblio]\label{Waterhouse} 
 W. C. Waterhouse. Automorphisms of the Bruhat ordering on
  Coxeter groups. {\it Bull. London Math. Soc.}, {\bf 21} : 243-248,
  1989.

\end{trivlist}

\end{document}